# Central limit theorems in linear structural error-in-variables models with explanatory variables in the domain of attraction of the normal law

## Yuliya V. Martsynyuk*

*School of Mathematics and Statistics, Carleton University,*
*1125 Colonel By Drive, Ottawa, ON K1S 5B6, Canada.*
*e-mail:* ymartsynyuk@yahoo.com

**Abstract:** Linear structural error-in-variables models with univariate observations are revisited for studying modified least squares estimators of the slope and intercept. New marginal central limit theorems (CLT's) are established for these estimators, assuming the existence of four moments for the measurement errors and that the explanatory variables are in the domain of attraction of the normal law. The latter condition for the explanatory variables is used the first time, and is so far the most general in this context. It is also optimal, or nearly optimal, for our CLT's. Moreover, due to the obtained CLT's being in Studentized and self-normalized forms to begin with, they are a priori nearly, or completely, data-based, and free of unknown parameters of the joint distribution of the error and explanatory variables. Consequently, they lead to a variety of readily available, or easily derivable, large-sample approximate confidence intervals (CI's) for the slope and intercept. In contrast, in related CLT's in the literature so far, the variances of the limiting normal distributions, in general, are complicated and depend on various, typically unknown, moments of the error and explanatory variables. Thus, the corresponding CI's for the slope and intercept in the literature, unlike those of the present paper, are available only under some additional model assumptions.

Received April 2007.

## Contents

*Research supported by a Carleton University Faculty of Graduate Studies and Research scholarship, an NSERC Canada Discovery Grant of M. Csörgő at Carleton University, and an NSERC Postdoctoral Fellowship of Yu.V. Martsynyuk at University of Ottawa.

## 1. Introduction

In the linear error-in-variables model (EIVM) of this paper we observe pairs $(y_i, x_i) \in \mathbb{R}^2$ according to

$$y_i = \beta \xi_i + \alpha + \delta_i, \quad x_i = \xi_i + \varepsilon_i, \tag{1}$$

where $\xi_i$ are unknown explanatory/latent variables, the real-valued slope $\beta$ and intercept $\alpha$ are to be estimated, and $\delta_i$ and $\varepsilon_i$ are unknown measurement error terms/variables, $1 \leq i \leq n$, $n \in \mathbb{N}$. EIVM (1) is also known as a measurement error model, or structural/functional relationship, or regression with errors in variables. It is a generalization of the simple linear regression of form $y_i = \beta \xi_i + \alpha + \delta_i$ in that in (1) it is also assumed that the two variables $\eta := \beta \xi + \alpha$ and $\xi$ are linearly related, and now not only $\eta$, but also $\xi$, are observed with respective measurement errors $\delta_i$ and $\varepsilon_i$.

In this paper we deal with the so-called structural EIVM (SEIVM), the model where the explanatory variables $\xi_i$ are assumed to be independent identically distributed (i.i.d.) random variables (r.v.'s) that are independent of the error terms (cf. **(C)** below). The case of (1) with $\alpha$ known to be zero is distinguished in the literature as the model without intercept. Convenient notations that are introduced in this paper (cf. (5), (6)) allow us to study both the no-intercept model and the model with unknown $\alpha$ simultaneously.

All our assumptions about (1) are summarized next:

**(A)** $\{(\delta, \varepsilon), (\delta_i, \varepsilon_i), i \geq 1\}$ is a sequence of i.i.d. random vectors with mean zero, $E \delta^4 < \infty$ and $E \varepsilon^4 < \infty$, and with a positive definite covariance matrix

$$\Gamma = \mathrm{Cov}(\delta, \varepsilon) = \begin{pmatrix} \lambda \theta & \mu \\ \mu & \theta \end{pmatrix}, \tag{2}$$

**(B)** $\{\xi, \xi_i, i \geq 1\}$ are i.i.d.r.v.'s in the domain of attraction of the normal law (DAN), i.e., there are constants $a_n$ and $b_n$, $b_n > 0$, for which $\left(\sum_{i=1}^n \xi_i - a_n\right) b_n^{-1} \xrightarrow{\mathcal{D}} N(0, 1)$, as $n \to \infty$,

**(C)** $(\delta, \varepsilon)$ and $\xi$ are independent.

**Remark 1.** Further to the definition of DAN in **(B)**, it is known that $a_n$ can be taken as $n E \xi$ and $b_n = n^{1/2} \ell_\xi(n)$, where $\ell_\xi(n)$ is a slowly varying function at infinity (i.e., $\ell_\xi(az)/\ell_\xi(z) \to 1$, as $z \to \infty$, for any $a > 0$), defined by the distribution of $\xi$. Moreover, $\ell_\xi(n) = \sqrt{\mathrm{Var}\,\xi} > 0$, if $\mathrm{Var}\,\xi < \infty$, and $\ell_\xi(n) \nearrow \infty$, as $n \to \infty$, if $\mathrm{Var}\,\xi = \infty$. Also, $\xi$ has moments of all orders less than 2, and the variance of $\xi$ is positive, but need not be finite.

Apart from making assumptions on the distribution of $(\xi, \delta, \varepsilon)$, to ensure identifiability of unknown parameters such as, for example, $\beta$ and $\alpha$ in model (1), it is common in the literature to make use of some side conditions in this regard, usually as conditions on the matrix $\Gamma$ of (2) in **(A)**. There are only a few frequently used identifiability assumptions, and here we deal with two of them that read as follows:

$$\text{Var}\,\delta = \lambda\theta \text{ and } \text{cov}(\delta, \varepsilon) = \mu \text{ are known, while } \text{Var}\,\varepsilon = \theta \text{ is unknown;} \quad (3)$$

$$\text{Var}\,\varepsilon = \theta \text{ and } \text{cov}(\delta, \varepsilon) = \mu \text{ are known, while } \text{Var}\,\delta = \lambda\theta \text{ is unknown.} \quad (4)$$

The identifiability assumption (4) is likely to be realistic in many applications. Neither (3) nor (4) correspond to orthogonal regression estimation in (1) that requires the assumption that the ratio $\text{Var}\,\delta/\text{Var}\,\varepsilon = \lambda$ is known and that $\text{cov}(\delta, \varepsilon) = 0$ (cf. Carroll and Ruppert (1996), Carroll *et at.* (2006), Cheng and Van Ness (1999), Fuller (1987) for some further discussions along these lines).

For further use throughout, for real-valued variables $\{u_i, 1 \le i \le n\}$ and $\{v_i, 1 \le i \le n\}$, we put

$$\overline{u} = \frac{1}{n}\sum_{i=1}^{n} u_i, \quad s_{i,uv} = (u_i - c\,\overline{u})(v_i - c\,\overline{v}) \quad \text{and} \quad S_{uv} = \frac{1}{n}\sum_{i=1}^{n} s_{i,uv}, \quad (5)$$

with constant

$$c = \begin{cases} 0 \text{ , if intercept } \alpha \text{ is known to be zero,} \\ 1 \text{ , if intercept } \alpha \text{ is unknown.} \end{cases} \quad (6)$$

Under (3), the modified least squares estimators for $\beta$ and $\alpha$ (cf. Cheng and Tsai (1995)) that have been used when it is assumed that $0 < \text{Var}\,\xi < \infty$ are

$$\widehat{\beta}_{1n} = \frac{S_{yy} - \lambda\theta}{S_{xy} - \mu} \quad \text{and} \quad \widehat{\alpha}_{1n} = \overline{y} - \overline{x}\,\widehat{\beta}_{1n}, \quad \text{if } S_{xy} - \mu \ne 0 \text{ and } S_{yy} - \lambda\theta > 0, \quad (7)$$

while those under (4) are

$$\widehat{\beta}_{2n} = \frac{S_{xy} - \mu}{S_{xx} - \theta} \quad \text{and} \quad \widehat{\alpha}_{2n} = \overline{y} - \overline{x}\,\widehat{\beta}_{2n}, \quad \text{if } S_{xx} - \theta > 0. \quad (8)$$

We note that the estimators in (7) and (8) coincide in form with the respective maximum likelihood estimators for $\beta$ and $\alpha$ that are derived and studied in the model (1) under the assumption that vector $(\xi, \delta, \varepsilon)$ is normally distributed.

In this paper, we revisit the estimators in (7) and (8), and prove several central limit theorems (CLT's) for each of them (cf. Theorems 1–3, and Corollary 1), under the distribution-free DAN assumption in **(B)** on the explanatory variables that, to the best of our knowledge, is the most general ever used so far in this context (cf. Remark 2). Moreover, according to Theorem 1, it is also optimal for the CLT's therein, and is nearly optimal for the CLT's in Theorems 2 and 3 (cf. Remark 3). As to the condition **(A)** on the error terms here, it seems to be the least restrictive that has been considered in the literature thus far.

Further to the special features of our CLT's in Theorems 1–3, all these CLT's are in Studentized or self-normalized forms to begin with and hence are automatically nearly, or completely, data-based. Namely, as compared to the CLT's for $\beta$ and $\alpha$ in the literature, Theorems 1–3 are a priori free of any unknown parameters of the distribution of $(\xi, \delta, \varepsilon)$ (cf. Remarks 3, 4), and, consequently, the corresponding large-sample approximate confidence intervals for the slope $\beta$ and intercept $\alpha$ are readily available, or easily derivable as in our main Theorem 4 (cf. the corresponding subsection right below Remark 8 for details).

Throughout Section 2, we pay a special attention to the SEIVM's (1) when $\text{Var}\,\xi = \infty$ (as allowed by **(B)** in view of Remark 1). Distinctive features of such models are obtained in Remarks 5–7 (cf. also Remark 9) that are seen to underpin our informal observation: the impact of the errors with finite variances in $x_i$ of (1) becomes automatically negligible as compared to that of the explanatory variables with $\text{Var}\,\xi = \infty$, and such SEIVM's are then close in spirit to, and behave like, the simple linear regression models $y_i = \beta x_i + \alpha + \delta_i$.

All the CLT's of this paper have strongly been inspired and influenced by recent advances in DAN via Studentization and self-normalization that are summarized in Csörgő, Szyszkowicz and Wang (2004). These developments prompted us to enrich the traditional two-moment space of the explanatory variables that has been used so far for CLT studies in SEIVM's (1) by allowing $\xi_i$ the first time to be simply in DAN with possibly infinite variance. For the use of DAN in some other regression models, we refer to Maller (1981) and Remark (iv) in Maller (1993).

## 2. Main results

We establish Studentized and self-normalized CLT's for each of the estimators in (7) and (8) in Theorems 1–3, and, together with the CLT's of Corollary 1 and the confidence intervals of Theorem 4, they constitute the main results of this paper.

**Theorem 1.** *Assume that the intercept $\alpha$ in (1) is known to be zero, $\{\xi, \xi_i, i \geq 1\}$ are i.i.d.r.v.'s, $E|\xi| < \infty$, and **(A)** and **(C)** hold true. Put*

$$U(j,n) = \begin{cases} S_{xy} - \mu \,, & \text{if } j = 1, \\ S_{xx} - \theta \,, & \text{if } j = 2, \end{cases} \quad u_i(j,n) = \begin{cases} (s_{i,yy} - \lambda\theta) - \beta(s_{i,xy} - \mu) \,, & \text{if } j = 1, \\ (s_{i,xy} - \mu) - \beta(s_{i,xx} - \theta) \,, & \text{if } j = 2. \end{cases}$$
(9)

*Then for $j = 1$ and $2$, **(B)** is equivalent to any one of the following CLT's: as $n \to \infty$,*

(a) $\sqrt{n}\,U(j,n)(\widehat{\beta}_{jn} - \beta)\Big(\sum_{i=1}^{n}(u_i(j,n) - \overline{u(j,n)})^2/(n-1)\Big)^{-1/2} \xrightarrow{\mathcal{D}} N(0,1)$,

(b) $n\,U(j,n)(\widehat{\beta}_{jn} - \beta)\Big(\sum_{i=1}^{n} u_i^2(j,n)\Big)^{-1/2} \xrightarrow{\mathcal{D}} N(0,1)$.

**Theorem 2.** *Assume* **(A)**–**(C)**. *Let*

$$\widetilde{u}_i(j,n) = \begin{cases} (s_{i,yy} - \lambda\theta) - \widehat{\beta}_{1n}(s_{i,xy} - \mu) & , \text{ if } j = 1, \\ (s_{i,xy} - \mu) - \widehat{\beta}_{2n}(s_{i,xx} - \theta) & , \text{ if } j = 2. \end{cases} \tag{10}$$

*Then, for $j = 1$ and $2$, the* CLT's *in the* (a) *and* (b) *parts of Theorem 1 hold true, and also,*

$$n\, U(j,n)(\widehat{\beta}_{jn} - \beta)\left(\sum_{i=1}^{n} \widetilde{u}_i^2(j,n)\right)^{-1/2} \xrightarrow{\mathcal{D}} N(0,1), \quad n \to \infty. \tag{11}$$

**Theorem 3.** *Assume* **(A)**–**(C)**. *For $j = 1$ and $2$, let*

$$\begin{aligned} v_i(j,n) &= (y_i - \alpha) - \beta x_i - \frac{\overline{x}}{U(j,n)}\, u_i(j,n), \\ \widetilde{v}_i(j,n) &= (y_i - \alpha) - \widehat{\beta}_{jn} x_i - \frac{\overline{x}}{U(j,n)}\, \widetilde{u}_i(j,n), \end{aligned} \tag{12}$$

*where $U(j,n)$, $u_i(j,n)$ and $\widetilde{u}_i(j,n)$ are as in* (9) *and* (10). *Then, for $j = 1$ and $2$, as $n \to \infty$,*

(a) $\sqrt{n}\,(\widehat{\alpha}_{jn} - \alpha)\left(\sum_{i=1}^{n}(v_i(j,n) - \overline{v(j,n)})^2/(n-1)\right)^{-1/2} \xrightarrow{\mathcal{D}} N(0,1),$

(b) $\sqrt{n}\,(\widehat{\alpha}_{jn} - \alpha)\left(\sum_{i=1}^{n}(\widetilde{v}_i(j,n) - \overline{\widetilde{v}(j,n)})^2/(n-1)\right)^{-1/2} \xrightarrow{\mathcal{D}} N(0,1).$

**Remark 2.** It follows from Cheng and Tsai (1995) that, under **(A)**, $0 < \text{Var}\,\xi < \infty$, **(C)**, independence of $\delta$ and $\varepsilon$, and $E\,\xi^4 < \infty$, $\widehat{\beta}_{jn}$ and $\widehat{\alpha}_{jn}$ are $\sqrt{n}$–asymptotically normal, for $j = 1$ and $2$. In Theorem 1.2.1 of Fuller (1987), among other things, the same conclusion is derived for $\widehat{\beta}_{2n}$ and $\widehat{\alpha}_{2n}$ under the condition that $(\xi, \delta, \varepsilon)$ is normally distributed with a positive definite diagonal covariance matrix. While Theorems 1–3 seem to present the first CLT's for SEIVM's (1) if $\text{Var}\,\xi = \infty$ (as allowed by **(B)** in view of Remark 1), they also imply the just mentioned CLT's for $\widehat{\beta}_{jn}$ and $\widehat{\alpha}_{jn}$ when $\text{Var}\,\xi < \infty$. Indeed, under the conditions of the specified CLT's in the literature, using the arguments from the proof of upcoming Corollary 1, expressions $U^{-2}(j,n)\sum_{i=1}^{n}(u_i(j,n) - \overline{u(j,n)})^2/(n-1)$ and $\sum_{i=1}^{n}(v_i(j,n) - \overline{v(j,n)})^2/(n-1)$ from the respective (a) parts of Theorems 1 and 3, $j = \overline{1,3}$, can be seen to converge in probability to positive constants that are the variances of the asymptotic normal distributions of the corresponding estimators obtained in Cheng and Tsai (1995) and Fuller (1987).

**Remark 3.** Due to their Studentized and self-normalized forms, the CLT's of Theorems 1–3 are invariant with respect to the distribution of $(\xi, \delta, \varepsilon)$ satisfying **(A)**–**(C)**, and are strikingly free of any unknown parameters of this distribution (depend only on the error moments that are assumed to be known in the identifiability assumptions (3) or (4)). In addition, (11) of Theorem 2 and the (b) part of Theorem 3 present completely data-based CLT's, while $\beta$ is the only unknown parameter appearing in the respective normalizers ($\sum_{i=1}^{n}(u_i(j,n) -$

$\overline{u(j,n)})^2/(n-1))^{-1/2}$, $\left(\sum_{i=1}^n u_i(j,n)^2\right)^{-1/2}$ and $(\sum_{i=1}^n (v_i(j,n)-\overline{v(j,n)})^2/(n-1))^{-1/2}$ of the CLT's of Theorem 1 that, according to Theorem 2, also holds true for the case of unknown $\alpha$, as well as of the CLT's of the (a) part of Theorem 3. Consequently, large-sample approximate confidence intervals for $\beta$ and $\alpha$ are readily available from (11) of Theorem 2 and the (b) part of Theorem 3, while those that follow from (a) and (b) of Theorem 1 are easily derivable. For the expressions for, and further discussions on, all these confidence intervals, we refer to the corresponding subsection, right below Remark 8. We also note that while **(B)** is optimal for the CLT's of Theorem 1 in the case of the no-intercept version of model (1), this condition is also optimal for the model (1) with an unknown intercept for the main terms in the expansions for the Studentized and self-normalized $\widehat{\beta}_{jn}$ and $\widehat{\alpha}_{jn}$ as in the CLT's of Theorems 2 and 3 (cf. respective proofs of Theorems 2 and 3, proof of Lemma 8, and conclusion of Lemma 7).

**Remark 4.** On account of their a priori Studentized and self-normalized forms, and also due to their respective features described in Remark 3, the CLT's of Theorem 1–3 appear to be also new when $\operatorname{Var}\xi < \infty$ (a special case of **(B)** in view of Remark 1). Indeed, as opposed to Theorems 1–3, in the CLT's for $\widehat{\beta}_{jn}$ and $\widehat{\alpha}_{jn}$ in Cheng and Tsai (1995), $j = 1$ and 2, that are proved under $\operatorname{Var}\xi < \infty$ (cf. Remark 2 for details), the expressions for the variances of the asymptotic normal distributions of $\widehat{\beta}_{jn}$ and $\widehat{\alpha}_{jn}$ are complicated and involve typically unknown, hard-to-estimate from data, moments of order $\leq 4$ of the error terms, in addition to the unknown parameters $\beta$, $E\xi$ and $\operatorname{Var}\xi$. Then, in order to be able to estimate the variances of the therein derived CLT's, it is additionally assumed that the errors $\delta$ and $\varepsilon$ are normally distributed. Consequently, the variances of the latter CLT's become simpler in form and contain only the unknown, but estimable $\beta$, $E\xi$, $\operatorname{Var}\xi$ and $\lambda\theta$ (or $\theta$). To handle similar difficulties with estimating the respective asymptotic variances of $\widehat{\beta}_{2n}$ and $\widehat{\alpha}_{2n}$ in Theorem 1.2.1 of Fuller (1987), the condition that $(\xi, \delta, \varepsilon)$ is normally distributed and has a positive definite diagonal covariance matrix is used. Consistent estimators for the respective asymptotic variances of $\widehat{\beta}_{jn}$ and $\widehat{\alpha}_{jn}$ that are proposed in the mentioned works are different from the expressions $U^{-2}(j,n)\sum_{i=1}^n \widetilde{u}_i^2(j,n)/n$ and $\sum_{i=1}^n (\widetilde{v}_i(j,n)-\overline{\widetilde{v}(j,n)})^2/(n-1)$, respectively taken from Theorems 2 and 3. When $\operatorname{Var}\xi < \infty$, due to (57) with $\eta_i(n) = u_i(j,n)$ and $\eta_i(n) = v_i'(j,n)$, (55), (58), (60), (61), (64), (68) and (69), the latter expressions appear to be the first consistent estimators for the just mentioned variances simply under **(A)**.

In Theorems 1–3, the rates of convergence to normality of $\widehat{\beta}_{jn}$ and $\widehat{\alpha}_{jn}$ are not apparent. For the sake of explicitly displaying these rates, we introduce the following direct consequence of Theorems 1–3.

**Corollary 1.** *If* **(A)**–**(C)** *are satisfied, then for $j = 1$ and 2, $\sqrt{n}\ell_\xi(n)(\widehat{\beta}_{jn} - \beta) \xrightarrow{\mathcal{D}} N(0, c_j)$ and $\sqrt{n}(\widehat{\alpha}_{jn} - \alpha) \xrightarrow{\mathcal{D}} N(0, d_j)$, as $n \to \infty$, where $\ell_\xi(n)$ is a typically unknown slowly varying function at infinity as in Remark 1 that converges to infinity when $\operatorname{Var}\xi = \infty$, and equals to a positive constant when $\operatorname{Var}\xi < \infty$, while $c_j$ and $d_j$ are positive constants.*

**Remark 5.** From Corollary 1, $\widehat{\beta}_{jn}$ are seen to be $\sqrt{n}\ell_{\xi}(n)$-asymptotically normal estimators of $\beta$. In this regard we note that when $\operatorname{Var}\xi = \infty$, the degree of precision of $\widehat{\beta}_{jn}$ increases as compared to the case $\operatorname{Var}\xi < \infty$. This effect rhymes well with our empirical expectation in that, intuitively, by letting $\xi_i$ in (1) to have an infinite deviation, we make them more dominant over the errors with finite variances. This, in turn, renders observations $y_i$ and $x_i$ to be more robust to noise (errors) and thus, more precise. As to the estimators $\widehat{\alpha}_{jn}$, according to Corollary 1, they are $\sqrt{n}$−asymptotically normal, regardless of whether $\operatorname{Var}\xi = \infty$, or $\operatorname{Var}\xi < \infty$.

**Remark 6.** We observe that if $\operatorname{Var}\xi = \infty$, then (3) and (4), as well as any other identifiability conditions, are unnecessary for constructing consistent estimators for $\beta$ and $\alpha$. For example, using similar arguments to those in (60) and (61), it can be shown that $S_{yy}/S_{xy}$ and $S_{xy}/S_{xx}$ are consistent estimators for $\beta$. The existence of these consistent estimators implies that $\beta$ is identifiable. The latter fact, when $(\delta, \varepsilon)$ has a normal distribution, can also be concluded from Reiersøl (1950), as accordingly, $\beta$ is identifiable if and only if $\xi$ is not normally distributed. As to consistent estimators for $\alpha$ under $\operatorname{Var}\xi = \infty$, one has $\overline{y} - \overline{x}\, S_{yy}/S_{xy}$ and $\overline{y} - \overline{x}\, S_{xy}/S_{xx}$.

**Remark 7.** Condition $\operatorname{Var}\xi = \infty$ can also be related to one of the frequently used identifiability assumptions for SEIVM's (1) that reads as follows:

$$\text{reliability ratio } k_{\xi} := \frac{E\xi^2 - c(E\xi)^2}{E\,\xi^2 - c(E\xi)^2 + \operatorname{Var}\varepsilon} \text{ is known,} \tag{13}$$
$$\text{under } \operatorname{Var}\xi < \infty \text{ and } E(\delta\varepsilon) = 0,$$

where $c$ is as in (6). In the case of $\operatorname{Var}\xi < \infty$, the coefficient $k_{\xi}$ plays a key role in the large sample theory of regression with errors in variables (1). In particular, $k_{\xi}$ adjusts the ordinary least squares estimator $S_{xy}/S_{xx}$ of the simple linear regression $y_i = \beta x_i + \alpha + \delta_i$ for consistency in (1) (that holds under **(A)** with $\mu = 0$ in (2), **(C)** and $0 < \operatorname{Var}\xi < \infty$) as follows: $k_{\xi}^{-1} S_{xy}/S_{xx}$. Now, defining $k_{\xi}$ of (13) to be 1 if $\operatorname{Var}\xi = \infty$, we have $k_{\xi}^{-1} S_{xy}/S_{xx} = S_{xy}/S_{xx}$, where the latter expression is one of the two proposed estimators for $\beta$ under $\operatorname{Var}\xi = \infty$ in Remark 6 that, in turn, coincides with the ordinary least squares estimator and does not require an adjustment via $k_{\xi}$ in (1) any more. The following view of the SEIVM (1) under $\operatorname{Var}\xi = \infty$ may shed light on this phenomenon. In $x_i$ of (1), the impact of the error terms with finite variances is negligible as compared to that of the explanatory variables with $\operatorname{Var}\xi = \infty$, and the model becomes close in spirit to, and behaves like, the simple linear regression $y_i = \beta x_i + \alpha + \delta_i$. We also note that if $\operatorname{Var}\xi < \infty$, then the reliability ratio $k_{\xi}$ usually has to be estimated from prior information for the sake of further use in inference in EIVM's. However, under $\operatorname{Var}\xi = \infty$, no estimation of $k_{\xi} := 1$ is necessary.

**Remark 8.** The present paper constitutes a part of Martsynyuk (2004), where, among other things, in the same context and spirit, the author studies weak/strong consistency and asymptotic normality of least squares estimators for the slope and intercept, as well as of methods of moments estimators for the error variances, under (3) and (4), and yet another identifiability condition that assumes that the matrix $\Gamma$ of (2) is known at least up to an unknown multiple $\theta = \operatorname{Var}\varepsilon$. As to the problem of proving new similarly featured CLT's for estimating $(\beta, \alpha)$ under the same model assumptions as those used in Theorem 1–3, namely under **(A)**–**(C)**, we again refer to Martsynyuk (2004). The author's Ph.D. thesis Martsynyuk (2005), among other things, also extends the just mentioned contributions of Martsynyuk (2004) regarding the SEIVM's to their traditional companions, the functional EIVM's (1), where the explanatory variables $\xi_i$ are assumed to be deterministic.

### Confidence intervals for slope $\beta$ and intercept $\alpha$

Abbreviations LSA and CI stand respectively for large-sample approximate and confidence interval, while $z_{\gamma/2}$ denotes the $100(1-\gamma/2)^{\text{th}}$ percentile of the standard normal distribution, $0 < \gamma < 1$.

In the SEIVM's (1) studied under (3) or (4), LSA CI's seem to be the only source of CI's for the slope $\beta$ and intercept $\alpha$. In order to work out completely data-based LSA CI's for $\beta$ and $\alpha$ from the corresponding CLT's in the literature, additional conditions of normality for the errors alone, or together with those on the explanatory variables have been used (cf. Remark 4). In contrast, on account of Theorems 1–3, computable LSA CI's for $\beta$ and $\alpha$ are now available under the very general distribution-free assumptions in **(A)**–**(C)** (cf. Remark 3), and also the first time under $\operatorname{Var}\xi = \infty$, as a special case of **(B)**. Thus, for $j = 1$ and 2, completely data-based CLT's of (11) of Theorem 2 and (b) of Theorem 3 imply readily available respective LSA $1 - \gamma$ CI's for $\beta$ and $\alpha$, $0 < \gamma < 1$, as follows:

$$\widehat{\beta}_{jn} \mp z_{\gamma/2} \frac{\left(\sum_{i=1}^{n} \widetilde{u}_i^2(j,n)\right)^{1/2}}{n\, U(j,n)}, \tag{14}$$

$$\widehat{\alpha}_{jn} \mp z_{\gamma/2} \frac{\left(\sum_{i=1}^{n} (\widetilde{v}_i(j,n) - \overline{\widetilde{v}(j,n)})^2\right)^{1/2}}{\sqrt{n(n-1)}}. \tag{15}$$

Moreover, under the usual assumption that $\operatorname{Var}\xi < \infty$, (14) and (15) are also new, being different from the aforementioned LSA CI's in the literature (cf. Remark 4).

Under each of the identifiability assumptions (3) and (4), in addition to the CI of (14), another two new LSA CI's for $\beta$ are also within reach. Namely, these are the LSA CI's for $\beta$ from the Studentized and self-normalized CLT's in (a) and (b) of Theorem 1 that, according to Theorem 2, also hold true for the SEIVM (1) with an unknown intercept $\alpha$. In these CLT's, $\beta$ is left unestimated in the normalizers $\left(\sum_{i=1}^{n}(u_i(j,n) - \overline{u(j,n)})^2/(n-1)\right)^{-1/2}$ and $\left(\sum_{i=1}^{n} u_i^2(j,n)/n\right)^{-1/2}$, respectively, $j = 1$ and 2, as opposed to the corresponding CLT's of (11) of Theorem 2 that are used to obtain (14). In Theorem 4, we obtain such CI's only under the identifiability assumption (3) (case $j = 1$). The case when that of (4) is being assumed can be handled similarly.

**Theorem 4.** *Assume* **(A)** *and* **(C)** *and that* $E|\xi|^{8/3} < \infty$*, which, in turn, implies* **(B)**. *Then, both for the no-intercept and unknown intercept versions of model* (1)*, from the* (a) *and* (b) *parts of Theorem 1 with* $j = 1$*, respective* LSA $1 - \gamma$ CI's *for* $\beta$ *are*

$$\left[ B_k^1(n, z_{\gamma/2}), B_k^2(n, z_{\gamma/2}) \right], \quad k = 1 \text{ and } 2, \tag{16}$$

*where for* $l = 1$ *and* 2,

$$B_k^l(n, z_{\gamma/2})$$
$$= \frac{f_k(n)(S_{xy} - \mu)^2 \widehat{\beta}_{1n} - z_{\gamma/2}^2 \sum_{i=1}^n (s_{i,yy} - g_k^{yy})(s_{i,xy} - g_k^{xy}) + (-1)^l \sqrt{\dfrac{D_k(n, z_{\gamma/2})}{4}}}{f_k(n)(S_{xy} - \mu)^2 - z_{\gamma/2}^2 \sum_{i=1}^n (s_{i,xy} - g_k^{xy})^2}, \tag{17}$$

*with*

$$f_1(n) = n(n-1), \ f_2(n) = n^2, \ g_1^{yy} = S_{yy}, \ g_2^{yy} = \lambda\theta, \ g_1^{xy} = S_{xy}, \ g_2^{xy} = \mu \tag{18}$$

*and*

$$D_k(n, z_{\gamma/2}) = 4z_{\gamma/2}^2 f_k(n)(S_{xy} - \mu)^2 \sum_{i=1}^n \left( (s_{i,yy} - g_k^{yy}) - \widehat{\beta}_{1n}(s_{i,xy} - g_k^{xy}) \right)^2$$

$$-4z_{\gamma/2}^4 \left( \sum_{i=1}^n (s_{i,yy} - g_k^{yy})^2 \sum_{i=1}^n (s_{i,xy} - g_k^{xy})^2 - \left( \sum_{i=1}^n (s_{i,yy} - g_k^{yy})(s_{i,xy} - g_k^{xy}) \right)^2 \right). \tag{19}$$

Individual and, in case of the CI's for $\beta$, also comparative performances of the obtained CI's in (14), (15), and (16), together with its corresponding analogues in terms of $\widehat{\beta}_{2n}$, are to be further investigated, and is a subject of the author's ongoing research.

**Remark 9.** Further to our discussions in Remark 7 on the reliability ratio $k_\xi$ defined in (13), we note that this coefficient $k_\xi$ has also played a key role in the literature so far in determining reasonably accurate LSA CI's for $\beta$ and $\alpha$ in SEIVM's (1). In particular, liabilities that some of LSA CI's for $\beta$ in SEIVM's (1) under $\mathrm{Var}\,\xi < \infty$ may suffer due to the so-called Gleser-Hwang effect (cf. Gleser and Hwang (1987)) are reasonably negligible in SEIVM's if the reliability ratio $k_\xi$ ($k_\xi < 1$) is far enough from zero (cf. Gleser (1987) for details). Though Gleser (1987) only deals with a specific LSA CI for $\beta$, assuming in this regard that the ratio of the uncorrelated variances is known, the just mentioned main conclusion is likely to be true for other available LSA CI's for $\beta$, and those for $\alpha$ in the SEIVM (1), when $\operatorname{Var} \xi < \infty$. In particular, it is desirable to rigorously support a reasonable belief that big enough $k_\xi$ will lead to accurate enough CI's in (14)–(16), as well as in the analogues of (16) in terms of $\widehat{\beta}_{2n}$, in the sense of having a negligible Gleser-Hwang effect. A common sense behind Gleser (1987) is that if $k_\xi = 0$, then $\operatorname{Var} \xi = 0$ and SEIVM (1) becomes degenerate, i.e., (1) reduces to $y_i = E\,\xi\,\beta + \alpha + \delta_i$ and $x_i = E\,\xi + \varepsilon_i$, where the explanatory variables do not vary any more, and thus it becomes impossible to fit a unique straight line through the data points. As opposed to the latter degenerate model, in SEIVM's (1) under $\operatorname{Var} \xi = \infty$, the explanatory variables are so well spread that they dominate the error terms in the sense that, according to our extended definition of $k_\xi$ in Remark 7, $k_\xi := 1$. Hence, it is only natural to conjecture that the Gleser-Hwang effect in regards of LSA CI's for $\beta$ and $\alpha$ disappears for such SEIVM's that in Remark 7 were seen to behave as if they were like the simple regression $y_i = x_i\beta + \alpha + \delta_i$.

## 3. Proofs of main results

### 3.1. Some auxiliary results

In this subsection we state some well-known results on DAN as Lemmas 1–4, give a simple alternative proof in the context of this paper for the first part of Lemma 4 on a characterization of DAN, and also establish a companion characterization as Lemma 5. Further developments in Section 3.2 leading to the proofs of the main results of Section 2 are based on these lemmas.

Hereafter, abbreviation WLLN stands for the Kolmogorov weak law of large numbers.

One of the several necessary and sufficient conditions for $\{Z, Z_i,\ i \geq 1\}$ to be in DAN is commonly associated with O'Brien (1980) as, e.g., in Giné, Götze and Mason (1997) (for more details see also Remark (iii) in Maller (1993), p.194), and it reads as follows.

**Lemma 1.** *For i.i.d.r.v.'s* $\{Z, Z_i,\ i \geq 1\}$, $Z \in \mathrm{DAN}$ *if and only if* $\max_{1 \leq i \leq n} Z_i^2 \big/ \sum_{i=1}^{n} Z_i^2 \overset{P}{\to} 0$, *as* $n \to \infty$.

The following result was rediscovered by Maller (1981), and is essentially a variation of Theorems 4 and 5 on pp. 143–144 in Gnedenko and Kolmogorov (1954).

**Lemma 2.** *Let* $\{Z, Z_i,\ i \geq 1\}$ *be i.i.d.r.v.'s in* DAN. *Then, as* $n \to \infty$, $\sum_{i=1}^{n}(Z_i - E\,Z)^2 b_n^{-2} \overset{P}{\to} 1$, *where* $b_n$ *is a positive sequence of numbers such that* $\sum_{i=1}^{n}(Z_i - E\,Z)b_n^{-1} \overset{\mathcal{D}}{\to} N(0,1)$.

The Giné, Götze and Mason (1997) fundamental characterization of DAN via Studentized or self-normalized partial sums can be stated as follows.

**Lemma 3.** *Let $\{Z, Z_i, i \geq 1\}$ be i.i.d.r.v.'s. Then, conditions $Z \in \mathrm{DAN}$ and $E\, Z = a$ are equivalent to any one of the following CLT's:* $\sqrt{n}(\overline{Z} - a) \left( \sum_{i=1}^{n}(Z_i - \overline{Z})^2/(n-1) \right)^{-1/2} \xrightarrow{\mathcal{D}} N(0, 1)$, *or* $\sum_{i=1}^{n}(Z_i - a) \left( \sum_{i=1}^{n}(Z_i - a)^2 \right)^{-1/2} \xrightarrow{\mathcal{D}} N(0, 1)$, *as $n \to \infty$.*

The first part of the following Lemma 4 that is due to Maller (1981) says that the DAN class of r.v.'s is closed under multiplication operation. Its second part amounts to a converse. Lemma 4 is applied in Maller (1981) to prove the asymptotic normality of the regression coefficient in a linear regression when the error variance is not necessarily finite. The proof of Lemma 4 in Maller (1981) is quite technical, and is based on checking the classical conditions of Theorem 2 on p. 128 in Gnedenko and Kolmogorov (1954) guaranteeing similar convergence in distribution to that in **(B)** for suitably chosen constants $b_n$. In the present context, we present a new, simpler and shorter proof of the first part of Lemma 4 under additionally assuming for $V \in \mathrm{DAN}$ therein that $E\, V^4 < \infty$. The latter assumption rhymes with our conditions $E\, \delta^4 < \infty$ and $E\, \varepsilon^4 < \infty$ in **(A)**, while assumption $U \in \mathrm{DAN}$ of the first part of Lemma 4 coincides with our condition **(B)** for $\xi$. The converse part of Lemma 4 is stated below without a proof.

**Lemma 4.** *Let $\{(U, V), (U_i, V_i), i \geq 1\}$ be i.i.d. random vectors, and $U$ and $V$ be independent. If $U \in \mathrm{DAN}$ and $V \in \mathrm{DAN}$, then $UV \in \mathrm{DAN}$. Conversely, if $E\, V^2 < \infty$ and $UV \in \mathrm{DAN}$, then $U \in \mathrm{DAN}$.*

*Proof.* . We prove here only the first part of Lemma 4, assuming additionally to $V \in \mathrm{DAN}$ that $E\, V^4 < \infty$.

If $E\, U^2 < \infty$, then, since $U$ and $V$ are independent and nondegenerate (both are in DAN), we have $\operatorname{Var} UV = E\, U^2 E\, V^2 - (E\, U)^2 (E\, V)^2 > 0$. Thus, on account of the CLT, $UV \in \mathrm{DAN}$.

Suppose now that $E\, U^2 = \infty$.

First, without loss of generality, we assume that $E\, V^2 = 1$ and prove the following key observation:

(i) $\sum_{i=1}^{n}(U_i - E\, U)^2 V_i^2 \Big/ \sum_{i=1}^{n}(U_i - E\, U)^2 \xrightarrow{P} 1, \quad n \to \infty$.

Since $U \in \mathrm{DAN}$, then $U - E\, U \in \mathrm{DAN}$ and, combining Lemma 3 and (3.7) of Giné, Götze and Mason (1997), one of the key results of that paper, we have

(ii) $E \left( \dfrac{|U_1 - E\, U|}{\left( \sum_{i=1}^{n}(U_i - E\, U)^2 \right)^{1/2}} \right)^4 = n^{-1} o(1)$.

For any $\varepsilon > 0$, on account of independence of $U$ and $V$, and (ii),

$$P\left( \left| \frac{\sum_{i=1}^{n}(U_i - E\, U)^2 V_i^2}{\sum_{i=1}^{n}(U_i - E\, U)^2} - 1 \right| > \varepsilon \right) \leq E\left( \sum_{i=1}^{n} \frac{(U_i - E\, U)^2}{\sum_{i=1}^{n}(U_i - E\, U)^2}(V_i^2 - 1) \right)^2 \varepsilon^{-2}$$

$$= n E(V_1^2 - 1)^2 E\left( (U_1 - E\, U)^2 \Big/ \sum_{i=1}^{n}(U_i - E\, U)^2 \right)^2 \varepsilon^{-2} = o(1),$$

i.e., we have (i).

Furthermore, without loss of generality, we can assume that $E\,U = 0$, since when $E\,U \neq 0$, it is easy to see that

(iii) $(U - E\,U)V \in \text{DAN}$   implies   $UV \in \text{DAN}$.

Indeed, if $(U - E\,U)V \in \text{DAN}$, then Lemma 2, (i) and the fact that $E\,U^2 = \infty$ yield, as $n \to \infty$,

(iv) $\dfrac{\sum_{i=1}^{n}(U_i - E\,U)V_i}{\sqrt{n}\,\ell(n)} \xrightarrow{\mathcal{D}} N(0,1),$

where the slowly varying function $\ell(n) \nearrow \infty$ is such that $(\sqrt{n}\ell(n))^{-1}\sum_{i=1}^{n}(U_i - E\,U) \xrightarrow{\mathcal{D}} N(0,1)$. Also, since $0 < \operatorname{Var} V < \infty$,

(v) $\dfrac{\sum_{i=1}^{n}(V_i - E\,V)E\,U}{\sqrt{n\operatorname{Var}V(E\,U)^2}} \xrightarrow{\mathcal{D}} N(0,1), \quad n \to \infty.$

Hence, on account of (iv) and (v), with $\ell(n)$ from (iv), as $n \to \infty$,

$$\frac{\sum_{i=1}^{n}(U_iV_i - E\,U \cdot E\,V)}{\sqrt{n}\,\ell(n)} = \frac{\sum_{i=1}^{n}(U_i - E\,U)V_i}{\sqrt{n}\,\ell(n)} + \frac{\sum_{i=1}^{n}(V_i - E\,V)E\,U}{\sqrt{n}\,\ell(n)}$$
$$\xrightarrow{\mathcal{D}} N(0,1),$$

i.e., $UV \in \text{DAN}$ via (iii) with $E\,U \neq 0$ and $E\,U^2 = \infty$.

Continuing the proof when $E\,U = 0$ and $E\,U^2 = \infty$, by Lemma 1, one needs to verify that

$$\max_{1 \leq i \leq n} U_i^2 V_i^2 \Big/ \sum_{i=1}^{n} U_i^2 V_i^2 = o_P(1), \quad \text{as } n \to \infty,$$

or, on account of (i), that

$$\max_{1 \leq i \leq n} U_i^2 V_i^2 \Big/ \sum_{i=1}^{n} U_i^2 = o_P(1), \quad \text{as } n \to \infty.$$

The latter, in turn, easily follows from Markov's inequality, independence of $U$ and $V$ and (ii), as follows: for every $\varepsilon > 0$,

$$P\left(\max_{1 \leq i \leq n} U_i^2 V_i^2 \Big/ \sum_{i=1}^{n} U_i^2 > \varepsilon\right) \leq nP\left(U_1^2 V_1^2 \Big/ \sum_{i=1}^{n} U_i^2 > \varepsilon\right)$$
$$\leq \varepsilon^{-2} n E\left(U_1^2 \Big/ \sum_{i=1}^{n} U_i^2\right)^2 E\,V_1^4 = o(1).$$

$\square$

**Remark 10.** Further to Lemma 4, from Remark on p.183 in Maller (1981) we learn that if $U \in \mathrm{DAN}$, $V \in \mathrm{DAN}$ and $E V^2 < \infty$, then $\sum_{i=1}^{n}(U_i - E U)V_i\left(\sqrt{n E V^2}\,\ell(n)\right)^{-1} \xrightarrow{\mathcal{D}} N(0,1)$ and $\sum_{i=1}^{n}(U_i - E U)\left(\sqrt{n}\,\ell(n)\right)^{-1} \xrightarrow{\mathcal{D}} N(0,1)$, $n \to \infty$, with the same slowly varying function at infinity $\ell(n)$, where, mutatis mutandis, $\ell(n)$ is as in Remark 1, i.e., accordingly featured for $\{U, U_i,\ i \geq 1\}$. A simple proof of this fact under additionally assuming that $E U^4 < \infty$ amounts to (iv), where $E V^2 = 1$.

The results of Lemma 4 will usually be coupled with those of the next one.

**Lemma 5.** *Let* $\{(U,V),(U_i,V_i),\ i \geq 1\}$ *be i.i.d. random vectors with* $E|UV| < \infty$.

  (a) *If* $U \in \mathrm{DAN}$ *and* $V \in \mathrm{DAN}$, *then* $U + V \in \mathrm{DAN}$, *provided that* $P(U + V = const) \neq 1$.

  (b) *Conversely, if* $U + V \in \mathrm{DAN}$, $U \in \mathrm{DAN}$ *and* $E U^2 < \infty$, *then* $V \in \mathrm{DAN}$, *provided that* $P(V = const) \neq 1$.

*Proof.* (a) If $E U^2 < \infty$, $E V^2 < \infty$ and $P(U + V = const) \neq 1$, then $0 < \mathrm{Var}\,(U + V) < \infty$ and, due to CLT, $U + V \in \mathrm{DAN}$.

Assume now that, say, $E U^2 = \infty$. Then, on account of Lemma 1, Lemma 2 for $\{U, U_i,\ i \geq 1\}$ with $b_n = \sqrt{n}\,\ell(n)$, where the slowly varying function $\ell(n) \nearrow \infty$ (cf. Remark 1 in terms of $\{U, U_i,\ i \geq 1\}$), and also the WLLN for $\{U_i V_i,\ i \geq 1\}$ with finite mean,

$$\frac{\max_{1 \leq i \leq n}(U_i + V_i)^2}{\sum_{i=1}^{n}(U_i + V_i)^2} \leq \frac{\max_{1 \leq i \leq n} 2(U_i^2 + V_i^2)}{\sum_{i=1}^{n} U_i^2\left(1 + 2\sum_{i=1}^{n} U_i V_i / \sum_{i=1}^{n} U_i^2\right) + \sum_{i=1}^{n} V_i^2}$$

$$\leq \frac{2\max_{1 \leq i \leq n} U_i^2}{\sum_{i=1}^{n} U_i^2(1 + o_P(1))} + \frac{2\max_{1 \leq i \leq n} V_i^2}{\sum_{i=1}^{n} V_i^2} = o_P(1).$$

Via Lemma 1, this proves that $U + V \in \mathrm{DAN}$.

  (b) Since $E U^2 < \infty$ and also $E|UV| < \infty$, then $E|(U + V)(-U)| < \infty$, and via applying the (a) part of Lemma 5 to $U + V$ and $-U$, we conclude that $V \in \mathrm{DAN}$. $\qquad \square$

### 3.2. Auxiliary results and proofs of main results

In the sequel, all vectors are row-vectors, and $\langle \cdot, \cdot \rangle$ stands for Euclidean inner product of two vectors. If $Z$ is a $d$−dimensional vector, then $Z^{(j)}$ is its $j^{\mathrm{th}}$ component, while $Z^{(k,k+l)} = (Z^{(k)}, Z^{(k+1)}, \cdots, Z^{(k+l)})$ is a subvector of $Z$ that has all the components of $Z$ starting with $Z^{(k)}$ and ending with $Z^{(k+l)}$, $1 \leq k \leq d-1$, $1 \leq l \leq k+l \leq d$, $d \geq 2$.

The proofs of the main results of Section 2 require some auxiliary results. First, in Lemmas 6 and 7, and Corollary 2, we will study Studentized and self-normalized partial sums that are based on i.i.d.r.v.'s $\{\langle \zeta_i,\ b \rangle,\ 1 \leq i \leq n\}$, where

$b \in \mathbb{R}^7$ is a nonzero vector of constants and

$$\zeta_i = \Big((\xi_i - c\,m)\delta_i,\, (\xi_i - c\,m)\varepsilon_i,\, \delta_i,\, \varepsilon_i,\, \delta_i\varepsilon_i - \mu,\, \delta_i^2 - \lambda\theta,\, \varepsilon_i^2 - \theta\Big), \quad 1 \leq i \leq n, \tag{20}$$

with constant $c$ as in (6) and $m := E\,\xi$. Such sums are the respective prototypes of the Studentized and self-normalized $\widehat{\beta}_{jn}$ as in (a) and (b) of Theorem 1 in the no-intercept version of model (1), $j = 1$ and 2. In (1) with an unknown intercept, the Studentized and self-normalized partial sums that are based on $\{\langle \eta_i(n), d \rangle, 1 \leq i \leq n\}$ play the same role, where

$$\eta_i(n) = (y_i - \alpha,\, x_i,\, s_{i,yy} - \lambda\theta,\, s_{i,xy} - \mu,\, s_{i,xx} - \theta), \quad 1 \leq i \leq n, \tag{21}$$

and a nonzero vector of constants $d \in \mathbb{R}^5$ is such that

$$d^{(1)}\beta + d^{(2)} = 0 \quad \text{and} \quad d^{(3)}\beta^2 + d^{(4)}\beta + d^{(5)} = 0. \tag{22}$$

The results in Lemma 8 and Corollary 3 for such partial sums will also be applied to derive (11) of Theorem 2 and Theorem 3 for $\widehat{\alpha}_{jn}$, $j = 1$ and 2. Moreover, the use of all the auxiliary results in this section can go beyond the immediate needs of this paper (cf. Remark 11). At the end of this subsection, we also prove Corollary 1 and obtain the CI's of Theorem 4.

Introduce vector

$$\zeta_0 = \Big((\xi - c\,m)\delta,\, (\xi - c\,m)\varepsilon,\, \delta,\, \varepsilon,\, \delta\varepsilon - \mu,\, \delta^2 - \lambda\theta,\, \varepsilon^2 - \theta\Big). \tag{23}$$

**Lemma 6.** *Assume* **(A)–(C)**. *When* $b^{(1)} = b^{(2)} = 0$, *assume additionally that*

$$\mathrm{Var}\,\langle \zeta_0, b \rangle > 0. \tag{24}$$

*Then, as* $n \to \infty$,

$$\frac{\sqrt{n}\,\overline{\langle \zeta, b \rangle}}{\left(\sum_{i=1}^n \langle \zeta_i - \overline{\zeta}, b \rangle^2/(n-1)\right)^{1/2}} \xrightarrow{\mathcal{D}} N(0,1) \;\; and \;\; \frac{n\,\overline{\langle \zeta, b \rangle}}{\left(\sum_{i=1}^n \langle \zeta_i, b \rangle^2\right)^{1/2}} \xrightarrow{\mathcal{D}} N(0,1). \tag{25}$$

*Proof.* By Lemma 3, (25) is valid if

$$\langle \zeta_0, b \rangle \in \mathrm{DAN} \tag{26}$$

and $E\,\langle \zeta_0, b \rangle = 0$. The latter equality follows from **(A)–(C)**. Hence, (26) needs to be shown.

If $b^{(1)} = b^{(2)} = 0$, then, since (24) holds true, we have $0 < \mathrm{Var}\,\langle \zeta_0, b \rangle < \infty$, and thus (26) as well.

Suppose now that $|b^{(1)}| + |b^{(2)}| > 0$. Then, on account of the first part of Lemma 4 and the fact that $\xi \in \mathrm{DAN}$ if and only if $\xi - c\,m \in \mathrm{DAN}$,

$$(\xi - c\,m)(b^{(1)}\delta + b^{(2)}\varepsilon) \in \mathrm{DAN}, \tag{27}$$

where, due to the fact that $\Gamma$ of (2) is positive definite,

$$\text{Var}(b^{(1)}\delta + b^{(2)}\varepsilon) > 0. \tag{28}$$

If $\text{Var}\Big(\langle\zeta_0, b\rangle - (\xi - c\,m)(b^{(1)}\delta + b^{(2)}\varepsilon)\Big) = 0$, then (27) implies (26). Otherwise, (26) is implied by part (a) of Lemma 5 applied to $(\xi - c\,m)(b^{(1)}\delta + b^{(2)}\varepsilon)$ and $\langle\zeta_0, b\rangle - (\xi - c\,m)(b^{(1)}\delta + b^{(2)}\varepsilon)$ that are now both from DAN. Two conditions of (a) of Lemma 5 are left to be verified. First, from the finiteness of $E\,\xi$ and the fourth error moments, and independence of $\xi$ and $(\delta, \varepsilon)$, it is seen that

$$
\begin{aligned}
&E\Big|(\xi - c\,m)(b^{(1)}\delta + b^{(2)}\varepsilon)\Big(\langle\zeta_0, b\rangle - (\xi - c\,m)(b^{(1)}\delta + b^{(2)}\varepsilon)\Big)\Big| \\
&= E\Big|(\xi - c\,m)(b^{(1)}\delta + b^{(2)}\varepsilon) \\
&\qquad \Big(b^{(3)}\delta + b^{(4)}\varepsilon + b^{(5)}(\delta\varepsilon - \mu) + b^{(6)}(\delta^2 - \lambda\theta) + b^{(7)}(\varepsilon^2 - \theta)\Big)\Big| < \infty.
\end{aligned}
\tag{29}
$$

If $\text{Var}\,\xi = \infty$, then the second assumption of the (a) part of Lemma 5, i.e., that $P(\langle\zeta_0, b\rangle = 0) \neq 1$, is automatically satisfied, since from $|b^{(1)}| + |b^{(2)}| > 0$, it follows that $\text{Var}\langle\zeta_0, b\rangle = \infty$.

It is left to be shown that $\text{Var}\langle\zeta_0, b\rangle > 0$ under assuming $\text{Var}\,\xi < \infty$ and $|b^{(1)}| + |b^{(2)}| > 0$. Consider the covariance matrix

$$
\text{Cov}\,\zeta_0 = \begin{pmatrix}
a_\xi\lambda\theta & a_\xi\mu & d_\xi\lambda\theta & d_\xi\mu & d_\xi m_{21} & d_\xi m_{30} & d_\xi m_{12} \\
a_\xi\mu & a_\xi\theta & d_\xi\mu & d_\xi\theta & d_\xi m_{12} & d_\xi m_{21} & d_\xi m_{03} \\
d_\xi\lambda\theta & d_\xi\mu & \lambda\theta & \mu & m_{21} & m_{30} & m_{12} \\
d_\xi\mu & d_\xi\theta & \mu & \theta & m_{12} & m_{21} & m_{03} \\
d_\xi m_{21} & d_\xi m_{12} & m_{21} & m_{12} & m_{22} - \mu^2 & m_{31} - \lambda\theta\mu & m_{13} - \theta\mu \\
d_\xi m_{30} & d_\xi m_{21} & m_{30} & m_{21} & m_{31} - \lambda\theta\mu & m_{40} - (\lambda\theta)^2 & m_{22} - \lambda\theta^2 \\
d_\xi m_{12} & d_\xi m_{03} & m_{12} & m_{03} & m_{13} - \theta\mu & m_{22} - \lambda\theta^2 & m_{04} - \theta^2
\end{pmatrix} \tag{30}
$$

where

$$a_\xi = E(\xi - c\,m)^2, \quad d_\xi = E(\xi - c\,m) \quad \text{and} \quad m_{ij} = E(\delta^i \varepsilon^j), \quad 0 \le i, j \le 4.$$

Let

$$
A = \begin{pmatrix}
\sqrt{a_\xi - d_\xi^2} & 0 & d_\xi & 0 \\
0 & \sqrt{a_\xi - d_\xi^2} & 0 & d_\xi \\
0 & 0 & 1 & 0 \\
0 & 0 & 0 & 1
\end{pmatrix} \tag{31}
$$

and $A^T$ denote the transpose of matrix $A$. By performing a straightforward multiplication, it can be verified that (30) is the product of three $7 \times 7$ block diagonal matrices, namely,

$$
\text{Cov}\,\zeta_0 = \begin{pmatrix} A & \mathrm{O} \\ \mathrm{O} & I_3 \end{pmatrix} \begin{pmatrix} \Gamma & \mathrm{O} \\ \mathrm{O} & \text{Cov}\,\zeta_0^{(3,7)} \end{pmatrix} \begin{pmatrix} A^T & \mathrm{O} \\ \mathrm{O} & I_3 \end{pmatrix} =: BCB^T, \tag{32}
$$

where O is a zero matrix of an appropriate size, $I_3$ is a $3 \times 3$ identity matrix and $\Gamma$ is as in (2). Since $\operatorname{Var}\xi > 0$ ($\xi \in \mathrm{DAN}$), $b^{(1,2)} \neq 0$ ($|b^{(1)}| + |b^{(2)}| > 0$) and $\Gamma$ is positive definite by **(A)**, then, on using (32), we conclude that

$$
\begin{aligned}
\operatorname{Var}\langle \zeta_0, b \rangle &= \langle b \operatorname{Cov} \zeta_0, b \rangle = \langle (bB)C, bB \rangle \\
&= \langle (bB)^{(1,2)}\Gamma, (bB)^{(1,2)} \rangle + \langle (bB)^{(3,7)}\operatorname{Cov}\zeta_0^{(3,7)}, (bB)^{(3,7)} \rangle \\
&\geq \langle (bB)^{(1,2)}\Gamma, (bB)^{(1,2)} \rangle = (a_\xi - d_\xi^2)\langle b^{(1,2)}\Gamma, b^{(1,2)} \rangle \\
&= \operatorname{Var}\xi \, \langle b^{(1,2)}\Gamma, b^{(1,2)} \rangle > 0.
\end{aligned}
$$

$\square$

**Corollary 2.** *Assume* **(A)**–**(C)***, and, when* $b^{(1)} = b^{(2)} = 0$*, also* (24)*. Then, as* $n \to \infty$*,*

$$
\overline{\langle \zeta, b \rangle} = \begin{cases}
\dfrac{O_P(1)}{\sqrt{n}}, & \text{if } \operatorname{Var}\xi < \infty \text{ and/or } b^{(1)} = b^{(2)} = 0, \\[2mm]
\dfrac{\ell_\xi(n)}{\sqrt{n}} O_P(1), & \text{if } \operatorname{Var}\xi = \infty \text{ and } |b^{(1)}| + |b^{(2)}| > 0,
\end{cases}
\tag{33}
$$

*where* $\ell_\xi(n)$ *is a slowly varying function at infinity as in Remark 1.*

*Proof.* If $\operatorname{Var}\xi < \infty$ and/or $b^{(1)} = b^{(2)} = 0$, then WLLN, (24) and positivity of $\operatorname{Var}\langle \zeta_0, b \rangle$ under $\operatorname{Var}\xi < \infty$ and $|b^{(1)}| + |b^{(2)}| > 0$, which is shown in the last paragraph of the proof of Lemma 6, result in

$$
\frac{1}{n-1}\sum_{i=1}^{n}\langle \zeta_i - \overline{\zeta}, b \rangle^2 \overset{P}{\to} \operatorname{Var}\langle \zeta_0, b \rangle > 0, \quad n \to \infty,
\tag{34}
$$

that combined with the first CLT in (25) yields $\overline{\langle \zeta, b \rangle} = O_P(1)/\sqrt{n}$.

If $\operatorname{Var}\xi = \infty$ and $|b^{(1)}| + |b^{(2)}| > 0$, then Remark 10 ($U = \xi$, $V = b^{(1)}\delta + b^{(2)}\varepsilon$), (28), Remark 1 and Lemma 2 imply

$$
\frac{\sum_{i=1}^{n}(\xi_i - m)^2(b^{(1)}\delta_i + b^{(2)}\varepsilon_i)^2}{n \ell_\xi^2(n)} \overset{P}{\to} \operatorname{Var}(b^{(1)}\delta + b^{(2)}\varepsilon) > 0, \quad n \to \infty,
$$

where $\ell_\xi(n) \nearrow \infty$. It can be shown that the latter convergence implies the convergence in probability of $\sum_{i=1}^{n}(\xi_i - c\,m)^2(b^{(1)}\delta_i + b^{(2)}\varepsilon_i)^2 / \big(n \ell_\xi^2(n)\big)$, $\sum_{i=1}^{n}\Big((\xi_i - c\,m)(b^{(1)}\delta_i + b^{(2)}\varepsilon_i) - \overline{(\xi - c\,m)(b^{(1)}\delta + b^{(2)}\varepsilon)}\Big)^2 / \big(n \ell_\xi^2(n)\big)$ and $\sum_{i=1}^{n}\langle \zeta_i - \overline{\zeta}, b \rangle^2 / \big(n \ell_\xi^2(n)\big)$, all to the same limit as above. From the first CLT in (25) and

$$
\frac{\sum_{i=1}^{n}\langle \zeta_i - \overline{\zeta}, b \rangle^2}{(n-1)\ell_\xi^2(n)} \overset{P}{\to} \operatorname{Var}(b^{(1)}\delta + b^{(2)}\varepsilon) > 0, \quad n \to \infty,
\tag{35}
$$

we have $\overline{\langle \zeta, b \rangle} = \ell_\xi(n) O_P(1)/\sqrt{n}$, where $\ell_\xi(n)/\sqrt{n} \to 0$, since $\ell_\xi(n)$ is slowly varying at infinity. $\square$

Next, we conclude that condition $\xi \in \mathrm{DAN}$ (cf. **(B)**) is optimal for any one of the CLT's in (25) when $|b^{(1)}| + |b^{(2)}| > 0$, i.e., when the Studentized and self-normalized partial sums in (25) are indeed based on $\{\xi_i,\ i \geq 1\}$, as opposed to the case $b^{(1)} = b^{(2)} = 0$, when they are error-based only.

**Lemma 7.** *Assume that for i.i.d.r.v.'s $\{\xi, \xi_i, i \geq 1\}$, $E|\xi| < \infty$. Let* **(A)** *and* **(C)** *hold true, and for vector $b \in \mathbb{R}^7$, $|b^{(1)}| + |b^{(2)}| > 0$. Then, $\xi \in \mathrm{DAN}$ is equivalent to any one of the following* CLT's*: $\sqrt{n}\ \overline{\langle \zeta, b \rangle} \Big( \sum_{i=1}^n \langle \zeta_i - \overline{\zeta},\ b \rangle^2 / (n-1) \Big)^{-1/2} \xrightarrow{\mathcal{D}} N(0,1)$, or $n\ \overline{\langle \zeta, b \rangle} \Big( \sum_{i=1}^n \langle \zeta_i, b \rangle^2 \Big)^{-1/2} \xrightarrow{\mathcal{D}} N(0,1)$, as $n \to \infty$.*

*Proof.* According to Lemma 3, any one of the above CLT's is equivalent to (26) and $E \langle \zeta_0, b \rangle = 0$. Due to the proof of Lemma 6, we only need to show that (26) implies that $\xi \in \mathrm{DAN}$.

If $\mathrm{Var}\Big( \langle \zeta_0, b \rangle - (\xi - c\,m)(b^{(1)}\delta + b^{(2)}\varepsilon) \Big) = 0$, then (26) immediately implies (27). Otherwise, (27) follows from the (b) part of Lemma 5 applied to $U + V = \langle \zeta_0, b \rangle$ and $U = \langle \zeta_0, b \rangle - (\xi - c\,m)(b^{(1)}\delta + b^{(2)}\varepsilon)$. Assumptions of this lemma are easily seen to be satisfied. We have $E|UV| < \infty$ on account of (29), and $U + V \in \mathrm{DAN}$ by (26). Clearly, $EU^2 = E\Big(b^{(3)}\delta + b^{(4)}\varepsilon + b^{(5)}(\delta\varepsilon - \mu) + b^{(6)}(\delta^2 - \lambda\theta) + b^{(7)}(\varepsilon^2 - \theta)\Big)^2 < \infty$, due to **(A)**. Finally, we are to verify that for $V = (\xi - c\,m)(b^{(1)}\delta + b^{(2)}\varepsilon)$, $P\Big(V = 0\Big) \neq 1$. The latter holds true since by Remark 1 and (28) that is due to having $|b^{(1)}| + |b^{(2)}| > 0$, $\mathrm{Var}\, V \geq \mathrm{Var}\, \xi \mathrm{Var}(b^{(1)}\delta + b^{(2)}\varepsilon) > 0$. Thus, we have checked conditions of the (b) part of Lemma 5 for obtaining (27). Now, from (27), **(A)** and the converse part of Lemma 4 we conclude that $\xi - c\,m \in \mathrm{DAN}$ and hence, $\xi \in \mathrm{DAN}$. $\qquad\square$

*Proof of Theorem* 1. The proof is due to Lemma 7 and the following representations:

$$\frac{\sqrt{n}U(j,n)(\widehat{\beta}_{jn} - \beta)}{\Big( \sum_{i=1}^n (u_i(j,n) - \overline{u(j,n)})^2 / (n-1) \Big)^{\frac{1}{2}}} = \frac{\sqrt{n}\ \overline{u(j,n)}}{\Big( \sum_{i=1}^n (u_i(j,n) - \overline{u(j,n)})^2 / (n-1) \Big)^{\frac{1}{2}}},$$

$$\frac{nU(j,n)(\widehat{\beta}_{jn} - \beta)}{\Big( \sum_{i=1}^n u_i^2(j,n) \Big)^{1/2}} = \frac{n\ \overline{u(j,n)}}{\Big( \sum_{i=1}^n u_i^2(j,n) \Big)^{1/2}} \qquad \text{and} \qquad u_i(j,n) = \langle \zeta_i, b_j \rangle,$$

$$\tag{36}$$

for $j = 1$ and 2, with $b_1 = (\beta, -\beta^2, 0, 0, -\beta, 1, 0)$ and $b_2 = (1, -\beta, 0, 0, 1, 0, -\beta)$ $(b_j^{(1,2)} \neq 0)$. $\qquad\square$

**Lemma 8.** *Assume* **(A)**–**(C)** *and* (22) *for a nonzero vector of constants $d \in \mathbb{R}^5$. Put*

$$e = (2\beta d^{(3)} + d^{(4)}, \beta d^{(4)} + 2d^{(5)}, d^{(1)}, d^{(2)}, d^{(4)}, d^{(3)}, d^{(5)}). \tag{37}$$

*When $e^{(1)} = e^{(2)} = 0$, assume additionally (24), with $e$ in place of $b$. Then, as $n \to \infty$,*

$$
\frac{\sqrt{n} \, \overline{\langle \eta(n), d \rangle}}{\left( \sum_{i=1}^{n} \langle \eta_i(n) - \overline{\eta(n)}, d \rangle^2 / (n-1) \right)^{1/2}} \xrightarrow{\mathcal{D}} N(0,1) \quad and
$$

$$
\frac{n \, \overline{\langle \eta(n), d \rangle}}{\left( \sum_{i=1}^{n} \langle \eta_i(n), d \rangle^2 \right)^{1/2}} \xrightarrow{\mathcal{D}} N(0,1).
$$

(38)

*Proof.* On account of (22),

$$
\begin{aligned}
\langle \eta_i(n), d \rangle &= d^{(1)}(\xi_i \beta + \delta_i) + d^{(2)}(\xi_i + \varepsilon_i) + d^{(3)}\Big(s_{i,\xi\xi}\beta^2 + 2s_{i,\xi\delta}\beta + (s_{i,\delta\delta} - \lambda\theta)\Big) \\
&\quad + d^{(4)}\Big(s_{i,\xi\xi}\beta + s_{i,\xi\delta} + s_{i,\xi\varepsilon}\beta + (s_{i,\delta\varepsilon} - \mu)\Big) + d^{(5)}\Big(s_{i,\xi\xi} + 2s_{i,\xi\varepsilon} + (s_{i,\varepsilon\varepsilon} - \theta)\Big) \\
&= (2\beta d^{(3)} + d^{(4)})(\xi_i - c\,\overline{\xi})(\delta_i - c\,\overline{\delta}) + (\beta d^{(4)} + 2d^{(5)})(\xi_i - c\,\overline{\xi})(\varepsilon_i - c\,\overline{\varepsilon}) \\
&\quad + d^{(1)}\delta_i + d^{(2)}\varepsilon_i + d^{(4)}\Big((\delta_i - c\,\overline{\delta})(\varepsilon_i - c\,\overline{\varepsilon}) - \mu\Big) \\
&\quad + d^{(3)}\Big((\delta_i - c\,\overline{\delta})^2 - \lambda\theta\Big) + d^{(5)}\Big((\varepsilon_i - c\,\overline{\varepsilon})^2 - \theta\Big) \\
&= \langle \zeta_i, e \rangle + c\,R_i(n),
\end{aligned}
$$

(39)

where vectors $\zeta_i$ and $e$ are as in (20) and (37), $c$ is from (6) and term $R_i(n)$ is

$$
\begin{aligned}
R_i(n) &= e^{(1)}\Big( -\overline{\delta}(\xi_i - m) + (m - \overline{\xi})(\delta_i - \overline{\delta}) \Big) \\
&\quad + e^{(2)}\Big( -\overline{\varepsilon}(\xi_i - m) + (m - \overline{\xi})(\varepsilon_i - \overline{\varepsilon}) \Big) + e^{(5)}\Big( -\overline{\varepsilon}\,\delta_i - \overline{\delta}\,\varepsilon_i + \overline{\delta}\,\overline{\varepsilon} \Big) \\
&\quad + e^{(6)}\Big( -2\overline{\delta}\,\delta_i + (\overline{\delta})^2 \Big) + e^{(7)}\Big( -2\overline{\varepsilon}\,\varepsilon_i + (\overline{\varepsilon})^2 \Big).
\end{aligned}
$$

(40)

If intercept $\alpha$ is known to be zero, i.e., $c = 0$, then the Studentized and self-normalized partial sums in (25) with vector $e$ in place of $b$, and respectively those in (38) coincide in view of (39). Therefore, the CLT's of Lemma 8 amount to those of Lemma 6.

Suppose now that $\alpha$ is unknown, i.e., $c = 1$. First, we will show that

$$
\begin{aligned}
\sqrt{n} \, \overline{R(n)} &= o_P(1) \text{ , if } \operatorname{Var} \xi < \infty \text{ and/or } e^{(1)} = e^{(2)} = 0, \\
\frac{\sqrt{n} \, \overline{R(n)}}{\ell_\xi(n)} &= o_P(1) \text{ , if } \operatorname{Var} \xi = \infty \text{ and } |e^{(1)}| + |e^{(2)}| > 0,
\end{aligned}
$$

(41)

with slowly varying function at infinity $\ell_\xi(n)$ as in Remark 1. For the summands in

$$
\sqrt{n} \, \overline{R(n)} = \sqrt{n}\Big( e^{(1)}\overline{\delta}(m - \overline{\xi}) + e^{(2)}\overline{\varepsilon}(m - \overline{\xi}) - e^{(5)}\overline{\delta}\,\overline{\varepsilon} - e^{(6)}(\overline{\delta})^2 - e^{(7)}(\overline{\varepsilon})^2 \Big),
$$

(42)

from the CLT for $\overline{\delta}$, **(B)** and Remark 1, as $n \to \infty$, we conclude

$$\frac{\sqrt{n}\,\overline{\delta}\,(\overline{\xi} - m)}{\ell_\xi(n)} = \frac{O_P(1)}{\sqrt{n}} \frac{\sqrt{n}(\overline{\xi} - m)}{\ell_\xi(n)} = \frac{O_P(1)}{\sqrt{n}} = o_P(1) \tag{43}$$

and, similarly,

$$\frac{\sqrt{n}\,\overline{\varepsilon}(m - \overline{\xi})}{\ell_\xi(n)} = o_P(1), \;\; \sqrt{n}\,\overline{\delta}\,\overline{\varepsilon} = o_P(1), \;\; \sqrt{n}(\overline{\delta})^2 = o_P(1), \;\; \sqrt{n}(\overline{\varepsilon})^2 = o_P(1). \tag{44}$$

This proves (41) that, combined with (34) and (35), yields

$$\frac{\sqrt{n}\,\overline{R(n)}}{\left(\sum_{i=1}^n \langle \zeta_i - \overline{\zeta}, e \rangle^2/(n-1)\right)^{1/2}} = o_P(1). \tag{45}$$

Now, in view of the first CLT in (25) of Lemma 6, (39) and (45), to complete the proof of the Studentized CLT in (38), it suffices to show that, as $n \to \infty$,

$$\frac{\sum_{i=1}^n \langle \eta_i(n) - \overline{\eta(n)}, d \rangle^2}{\sum_{i=1}^n \langle \zeta_i - \overline{\zeta}, e \rangle^2} \xrightarrow{P} 1, \tag{46}$$

where

$$\langle \zeta_i - \overline{\zeta}, e \rangle = e^{(1)}\Big((\xi_i - m)\delta_i - \overline{(\xi - m)\delta}\Big) + e^{(2)}\Big((\xi_i - m)\varepsilon_i - \overline{(\xi - m)\varepsilon}\Big)$$
$$+ e^{(3)}(\delta_i - \overline{\delta}) + e^{(4)}(\varepsilon_i - \overline{\varepsilon}) + e^{(5)}(\delta_i\varepsilon_i - \overline{\delta\varepsilon}) + e^{(6)}(\delta_i^2 - \overline{\delta^2}) + e^{(7)}(\varepsilon_i^2 - \overline{\varepsilon^2}), \tag{47}$$

and, as a consequence of (39),

$$\langle \eta_i(n) - \overline{\eta(n)}, d \rangle = e^{(1)}\Big((\xi_i - \overline{\xi})(\delta_i - \overline{\delta}) - (\overline{\xi\delta} - \overline{\xi}\,\overline{\delta})\Big)$$
$$+ e^{(2)}\Big((\xi_i - \overline{\xi})(\varepsilon_i - \overline{\varepsilon}) - (\overline{\xi\varepsilon} - \overline{\xi}\,\overline{\varepsilon})\Big)$$
$$+ e^{(3)}(\delta_i - \overline{\delta}) + e^{(4)}(\varepsilon_i - \overline{\varepsilon}) + e^{(5)}\Big((\delta_i - \overline{\delta})(\varepsilon_i - \overline{\varepsilon}) - (\overline{\delta\varepsilon} - \overline{\delta}\,\overline{\varepsilon})\Big)$$
$$+ e^{(6)}\Big((\delta_i - \overline{\delta})^2 - (\overline{\delta^2} - (\overline{\delta})^2)\Big) + e^{(7)}\Big((\varepsilon_i - \overline{\varepsilon})^2 - (\overline{\varepsilon^2} - (\overline{\varepsilon})^2)\Big). \tag{48}$$

By (34), (35) and the Cauchy-Schwarz inequality, to show (46), it suffices to prove that

$$\frac{1}{n}\sum_{i=1}^n \Big(\langle \eta_i(n) - \overline{\eta(n)}, d \rangle - \langle \zeta_i - \overline{\zeta}, e \rangle\Big)^2 = o_P(1), \quad n \to \infty. \tag{49}$$

On using the Cauchy-Schwarz inequality again, (49) follows from the following statements for the corresponding summands in (47) and (48): as $n \to \infty$,

$$\frac{1}{n}\sum_{i=1}^n \left(e^{(1)}\Big((\xi_i - \overline{\xi})(\delta_i - \overline{\delta}) - (\overline{\xi\delta} - \overline{\xi}\,\overline{\delta})\Big) - e^{(1)}\Big((\xi_i - m)\delta_i - \overline{(\xi - m)\delta}\Big)\right)^2 = o_P(1), \tag{50}$$

$$\frac{1}{n}\sum_{i=1}^{n}\left(e^{(2)}\Big((\xi_i-\overline{\xi})(\varepsilon_i-\overline{\varepsilon})-(\overline{\xi\varepsilon}-\overline{\xi}\,\overline{\varepsilon})\Big)-e^{(2)}\Big((\xi_i-m)\varepsilon_i-\overline{(\xi-m)\varepsilon}\Big)\right)^2=o_P(1),$$
(51)

$$\frac{1}{n}\sum_{i=1}^{n}\left(e^{(5)}\Big((\delta_i-\overline{\delta})(\varepsilon_i-\overline{\varepsilon})-(\overline{\delta\varepsilon}-\overline{\delta}\,\overline{\varepsilon})\Big)-e^{(5)}(\delta_i\varepsilon_i-\overline{\delta\varepsilon})\right)^2=o_P(1),$$
(52)

$$\frac{1}{n}\sum_{i=1}^{n}\left(e^{(6)}\Big((\delta_i-\overline{\delta})^2-(\overline{\delta^2}-(\overline{\delta})^2)\Big)-e^{(6)}(\delta_i^2-\overline{\delta^2})\right)^2=o_P(1)$$
(53)

and

$$\frac{1}{n}\sum_{i=1}^{n}\left(e^{(7)}\Big((\varepsilon_i-\overline{\varepsilon})^2-(\overline{\varepsilon^2}-(\overline{\varepsilon})^2)\Big)-e^{(7)}(\varepsilon_i^2-\overline{\varepsilon^2})\right)^2=o_P(1).$$
(54)

Similarly, (50) holds true on account of having

$$\frac{1}{n}\sum_{i=1}^{n}\Big((\xi_i-\overline{\xi})(\delta_i-\overline{\delta})-(\xi_i-m)\delta_i\Big)^2$$

$$=\frac{1}{n}\sum_{i=1}^{n}\Big(-(\xi_i-m)\overline{\delta}+(m-\overline{\xi})(\delta_i-\overline{\delta})\Big)^2$$

$$\leq\frac{2}{n}\sum_{i=1}^{n}(\xi_i-m)^2(\overline{\delta})^2+\frac{2}{n}\sum_{i=1}^{n}(\delta_i-\overline{\delta})^2(m-\overline{\xi})^2$$

$$=O_P(1)\,\frac{1}{n^2}\sum_{i=1}^{n}(\xi_i-m)^2+o_P(1)=o_P(1)$$

and

$$\frac{1}{n}\sum_{i=1}^{n}\left(\overline{(\xi-m)\delta}-(\overline{\xi\delta}-\overline{\xi}\,\overline{\delta})\right)^2=\frac{n(m-\overline{\xi})^2(\overline{\delta})^2}{n}=o_P(1),$$

where we have applied WLLN, the CLT for $\overline{\delta}$ and the Marcinkiewicz-Zygmund law of large numbers for $(\xi_i-m)^2$, where $E|(\xi_i-m)^2|^{1/2}<\infty$. (51) is obtained in the same manner. All (52)–(54) are handled similarly, and easily result from the Cauchy-Schwarz inequality and the WLLN under **(A)**. This completes the proof of (46), and hence also that of the first CLT in (38). The latter CLT implies that, as $n\to\infty$,

$$\frac{n\overline{(\langle\eta(n),d\rangle)^2}}{\sum_{i=1}^{n}\langle\eta_i(n)-\overline{\eta(n)},d\rangle^2}=\frac{O_P(1)}{n}=o_P(1),$$

which combined with the Cauchy-Schwarz inequality proves that

$$\frac{\sum_{i=1}^{n}\langle\eta_i(n),d\rangle^2}{\sum_{i=1}^{n}\langle\eta_i(n)-\overline{\eta(n)},d\rangle^2}\xrightarrow{P}1,\quad n\to\infty.$$
(55)

The Studentized CLT in (38) and (55) lead to the second, self-normalized CLT in (38). $\qquad\square$

**Corollary 3.** *Let all the assumptions of Lemma 8 be satisfied. Then, as $n \to \infty$,*

$$\overline{\langle \eta(n), d \rangle} = \begin{cases} \dfrac{\ell_\xi(n) O_P(1)}{\sqrt{n}} & , \ if \ |e^{(1)}| + |e^{(2)}| > 0, \\ \dfrac{O_P(1)}{\sqrt{n}} & , \ if \ e^{(1)} = e^{(2)} = 0, \end{cases} \tag{56}$$

*with slowly varying function at infinity $\ell_\xi(n)$ as in Remark 1, and vector $e$ of (37).*

*Proof.* From (34), (35) and (46), as $n \to \infty$,

$$\frac{\sum_{i=1}^{n} \langle \eta_i(n) - \overline{\eta(n)}, d \rangle^2}{(n-1)\ell_\xi^2(n)} \xrightarrow{P} const > 0, \ \text{ if } \mathrm{Var}\, \xi = \infty \text{ and } |e^{(1)}| + |e^{(2)}| > 0,$$

$$\frac{\sum_{i=1}^{n} \langle \eta_i(n) - \overline{\eta(n)}, d \rangle^2}{n-1} \xrightarrow{P} const > 0, \text{ otherwise.}$$

$$\tag{57}$$

The first CLT in (38) and (57) result in (56). □

**Remark 11.** Lemmas 7, 8 and Corollaries 2, 3 are rather versatile and, apart from the needs of this paper, can also be applied to establish Studentized and self-normalized marginal CLT's for other estimators that are appropriately based on the vector $(\overline{y}, \overline{x}, S_{yy}, S_{xy}, S_{xx})$ in the context of the SEIVM (1) (cf., e.g., such CLT's for the weighted least squares estimators for $\beta$ and $\alpha$, and for methods of moments estimators for the error variances $\lambda\theta$ and $\theta$ proved in Martsynyuk (2004)).

*Proof of* (a) *and* (b) *of Theorem 1 under the conditions of Theorem 2.* In view of Theorem 1, we only need to argue that (a) and (b) of Theorem 1 also hold true in the model (1) with an unknown intercept, provided **(A)**–**(C)** are assumed. This is a consequence of Lemma 8 and the representations in (36) for $j = 1$ and 2, where now $u_i(j, n) = \langle \eta_i(n), d_j \rangle$, with $d_1 = (0, 0, 1, -\beta, 0)$ and $d_2 = (0, 0, 0, 1, -\beta)$ that satisfy (22), and with vectors $e_1$ and $e_2$ corresponding to $e$ as in (37) that are equal to $b_1$ and $b_2$ as specified in the proof of Theorem 1. □

*Proof of* (11) *of Theorem 2.* In view of (b) of Theorem 1 that also holds under the conditions of Theorem 2, it suffices to show that, for $j = 1$ and 2,

$$\frac{\sum_{i=1}^{n} \widetilde{u}_i^2(j, n)}{\sum_{i=1}^{n} u_i^2(j, n)} \xrightarrow{P} 1, \quad n \to \infty. \tag{58}$$

We first note that, on using Remark 1 and Lemma 2,

$$\frac{\sum_{i=1}^{n} \xi_i^4}{n^2 \ell_\xi^4(n)} \leq 8 \left( \frac{\sum_{i=1}^{n} (\xi_i - m)^4}{n^2 \ell_\xi^4(n)} + \frac{n m^4}{n^2 \ell_\xi^4(n)} \right)$$

$$= O_P(1) \frac{\sum_{i=1}^{n} (\xi_i - m)^4}{\left( \sum_{i=1}^{n} (\xi_i - m)^2 \right)^2} + o_P(1) = o_P(1), \tag{59}$$

with $m = E\xi$ and $\ell_\xi(n)$ of Remark 1, since, due to (3.7) in Giné, Götze and Mason (1997) and identical distribution of the r.v.'s $(\xi_i - m)^4 \left( \sum_{i=1}^n (\xi_i - m)^2 \right)^{-2}$, $1 \le i \le n$, for any $\varepsilon > 0$,

$$P\left( \frac{\sum_{i=1}^n (\xi_i - m)^4}{(\sum_{i=1}^n (\xi_i - m)^2)^2} > \varepsilon \right) \le \varepsilon^{-1} n E \left( \frac{(\xi_1 - m)^4}{(\sum_{i=1}^n (\xi_i - m)^2)^2} \right) \to 0, \quad n \to \infty.$$

Also, via Remark 1, Lemma 2 and WLLN, as $n \to \infty$,

$$\frac{S_{xx} - \theta}{\ell_\xi^2(n)} = \frac{S_{\xi\xi}}{\ell_\xi^2(n)} + \frac{2S_{\xi\varepsilon}}{\ell_\xi^2(n)} + \frac{S_{\varepsilon\varepsilon} - \theta}{\ell_\xi^2(n)} = \frac{S_{\xi\xi}}{\ell_\xi^2(n)} + o_P(1)$$

$$= \frac{\overline{(\xi - m)^2} + 2\overline{(\xi - m)}(m - c\overline{\xi}) + (m - c\overline{\xi})^2}{\ell_\xi^2(n)} + o_P(1)$$

$$\xrightarrow{P} \begin{cases} \dfrac{E\xi^2 - c\,m^2}{\operatorname{Var}\xi}, & \text{if } \operatorname{Var}\xi < \infty, \\ 1, & \text{if } \operatorname{Var}\xi = \infty. \end{cases} \tag{60}$$

Similarly,

$$\frac{S_{xy} - \mu}{\beta\ell_\xi^2(n)} \xrightarrow{P} \begin{cases} \dfrac{E\xi^2 - c\,m^2}{\operatorname{Var}\xi}, & \text{if } \operatorname{Var}\xi < \infty, \\ 1, & \text{if } \operatorname{Var}\xi = \infty. \end{cases} \tag{61}$$

By (36), Corollary 3 (case of $|e^{(1)}| + |e^{(2)}| > 0$), (59), (60) and WLLN,

$$\frac{\sum_{i=1}^n (\widetilde{u}_i(2, n) - u_i(2, n))^2}{n\ell_\xi^2(n)} = \frac{\sum_{i=1}^n (s_{i,xx} - \theta)^2 (\beta - \widehat{\beta}_{2n})^2}{n\ell_\xi^2(n)}$$

$$= \left( \frac{\overline{u(2, n)}}{S_{xx} - \theta} \right)^2 \frac{\sum_{i=1}^n (s_{i,xx} - \theta)^2}{n\ell_\xi^2(n)} = O_P(1) \frac{\sum_{i=1}^n (s_{i,xx} - \theta)^2}{n^2 \ell_\xi^4(n)}$$

$$\le O_P(1) \frac{\sum_{i=1}^n s_{i,xx}^2}{n^2 \ell_\xi^4(n)} + o_P(1) \le O_P(1) \frac{\sum_{i=1}^n (x_i^4 + (\overline{x})^4)}{n^2 \ell_\xi^4(n)} + o_P(1)$$

$$\le O_P(1) \frac{\sum_{i=1}^n (\xi_i^4 + \varepsilon_i^4)}{n^2 \ell_\xi^4(n)} + \frac{n O_P(1) o_P(1)}{n^2 \ell_\xi^4(n)} + o_P(1) = o_P(1). \tag{62}$$

For $j = 2$, (58) follows from (36) with $u_i(2, n) = \langle \eta_i(n), d_2 \rangle$, $d_2 = (0, 0, 0, 1, -\beta)$, (55), (57) (case of $|e^{(1)}| + |e^{(2)}| > 0$), (62) and the Cauchy-Schwarz inequality.

The proof of (58) for $j = 1$ is similar to the above lines. $\qquad \square$

*Proof of Corollary 1 for $\widehat{\beta}_{jn}$.* The proof follows from (36) for $j = 1$ and 2, with $u_i(j, n) = \langle \eta_i(n), d_j \rangle$, where $d_1 = (0, 0, 1, -\beta, 0)$ and $d_2 = (0, 0, 0, 1, -\beta)$, and from Corollary 3 (case of $|e^{(1)}| + |e^{(2)}| > 0$), (60) and (61). $\qquad \square$

*Proof of the* (a) *part of Theorem* 3. Below we consider the case of $\widehat{\alpha}_{2n}$ only, as the respective CLT for $\widehat{\alpha}_{1n}$ can be proved in a similar way.

Suppose first that $\operatorname{Var}\xi = M < \infty$. Then, by Corollary 3 for $\overline{u(2,n)}$, (60), Remark 1 and WLLN, as $n \to \infty$,

$$
\begin{aligned}
\sqrt{n}(\widehat{\alpha}_{2n} - \alpha) &= \sqrt{n}\Big(\overline{(y - x\beta - \alpha)} - \overline{x}(\widehat{\beta}_{2n} - \beta)\Big) \\
&= \sqrt{n}\left(\overline{(y - x\beta - \alpha)} - \frac{m}{M}\,\overline{u(2,n)} + \left(\frac{m}{M} - \frac{\overline{x}}{S_{xx} - \theta}\right)\overline{u(2,n)}\right) \\
&= \sqrt{n}\,\overline{v'(2,n)} + o_P(1),
\end{aligned}
\tag{63}
$$

where $m = E\,\xi$ and

$$
v'_i(2,n) = y_i - \alpha - \beta x_i - \frac{m}{M}u_i(2,n).
$$

For $\sqrt{n}\,\overline{v'(2,n)}\Big(\sum_{i=1}^{n}(v'_i(2,n) - \overline{v'(2,n)})^2/(n-1)\Big)^{-1/2}$, (38) of Lemma 8 holds true, since the respective vector $d$ satisfies (22), and $e$ of (37) is such that $|e^{(1)}| + |e^{(2)}| > 0$ if $m \neq 0$, and if $m = 0$, $e^{(1)} = e^{(2)} = 0$ with (24), i.e., with inequality $\operatorname{Var}(\delta - \beta\varepsilon) > 0$ that is satisfied on account of (**A**). Combining this CLT, (57) with $\eta_i(n) = v'_i(2,n)$ and (63), we conclude that

$$
\sqrt{n}(\widehat{\alpha}_{2n} - \alpha)\left(\sum_{i=1}^{n}(v'_i(2,n) - \overline{v'(2,n)})^2/(n-1)\right)^{-1/2} \xrightarrow{\mathcal{D}} N(0,1), \quad n \to \infty.
$$

Hence, to complete the proof when $\operatorname{Var}\xi < \infty$, the following convergence has to be shown: as $n \to \infty$,

$$
\begin{aligned}
&\frac{\sum_{i=1}^{n}(v_i(2,n) - \overline{v(2,n)})^2}{\sum_{i=1}^{n}(v'_i(2,n) - \overline{v'(2,n)})^2} \\
&= \frac{\sum_{i=1}^{n}\Big((y_i - \overline{y}) - \beta(x_i - \overline{x}) - \frac{\overline{x}}{S_{xx}-\theta}(u_i(2,n) - \overline{u(2,n)})\Big)^2}{\sum_{i=1}^{n}\Big((y_i - \overline{y}) - \beta(x_i - \overline{x}) - \frac{m}{M}(u_i(2,n) - \overline{u(2,n)})\Big)^2} \xrightarrow{P} 1.
\end{aligned}
\tag{64}
$$

On observing that WLLN, (60) and (57) with $\eta_i(n) = u_i(2,n)$ imply

$$
\left(\frac{m}{M} - \frac{\overline{x}}{S_{xx} - \theta}\right)^2 \frac{1}{n-1}\sum_{i=1}^{n}(u_i(2,n) - \overline{u(2,n)})^2 = o_P(1), \quad n \to \infty,
$$

the proof of (64) follows from the Cauchy-Schwarz inequality and (57) with $\eta_i(n) = v'_i(2,n)$.

Assume now that $\operatorname{Var}\xi = \infty$. As $n \to \infty$, we are to prove convergence to $N(0,1)$ of

$$
\frac{\sqrt{n}(\widehat{\alpha}_{2n} - \alpha)}{\Big(\sum_{i=1}^{n}(v_i(2,n) - \overline{v(2,n)})^2/(n-1)\Big)^{1/2}}
$$

$$= \frac{\sqrt{n}\,\overline{v'(2,n)} - \sqrt{n}\,\overline{x}(\widehat{\beta}_{2n} - \beta)}{\left(\sum_{i=1}^{n}(v_i'(2,n) - \overline{v'(2,n)})^2/(n-1)\right)^{1/2}} \left(\frac{\sum_{i=1}^{n}(v_i'(2,n) - \overline{v'(2,n)})^2}{\sum_{i=1}^{n}(v_i(2,n) - \overline{v(2,n)})^2}\right)^{1/2},$$

(65)

where in this case

$$v_i'(2,n) = y_i - \alpha - \beta x_i.$$

As $n \to \infty$, by Corollary 1 for $\widehat{\beta}_{2n}$ and WLLN,

$$\sqrt{n}\,\overline{x}(\widehat{\beta}_{2n} - \beta) = \frac{O_P(1)}{\ell_\xi(n)} = o_P(1), \tag{66}$$

with $\ell_\xi(n)$ as in Remark 1, while **(A)** and WLLN give

$$\frac{1}{n-1}\sum_{i=1}^{n}(v_i'(2,n) - \overline{v'(2,n)})^2 = \frac{1}{n-1}\sum_{i=1}^{n}\left((y_i - \overline{y}) - \beta(x_i - \overline{x})\right)^2$$

$$= \frac{1}{n-1}\sum_{i=1}^{n}\left((\delta_i - \overline{\delta}) - \beta(\varepsilon_i - \overline{\varepsilon})\right)^2 \xrightarrow{P} \mathrm{Var}(\delta - \beta\varepsilon) > 0. \tag{67}$$

By (60), the WLLN for $\overline{x}$ and (57) for $\eta_i(n) = u_i(2,n)$ (case of $\mathrm{Var}\,\xi = \infty$, $|e^{(1)}| + |e^{(2)}| > 0$),

$$\left(\frac{\overline{x}}{S_{xx} - \theta}\right)^2 \frac{1}{n-1}\sum_{i=1}^{n}(u_i(2,n) - \overline{u(2,n)})^2 = \frac{O_P(1)}{\ell_\xi^2(n)} = o_P(1)$$

and hence, via the Cauchy-Schwarz inequality and (67), as $n \to \infty$,

$$\frac{\sum_{i=1}^{n}(v_i(2,n) - \overline{v(2,n)})^2}{\sum_{i=1}^{n}(v_i'(2,n) - \overline{v'(2,n)})^2}$$

$$= \frac{\sum_{i=1}^{n}\left((y_i - \overline{y}) - \beta(x_i - \overline{x}) - \frac{\overline{x}}{S_{xx} - \theta}(u_i(2,n) - \overline{u(2,n)})\right)^2}{\sum_{i=1}^{n}\left((y_i - \overline{y}) - \beta(x_i - \overline{x})\right)^2} \xrightarrow{P} 1. \tag{68}$$

Finally, combining (65)–(68) and (38) of Lemma 8 for $\sqrt{n}\,\overline{v'(2,n)}\left(\sum_{i=1}^{n}(v_i'(2,n) - \overline{v'(2,n)})^2/(n-1)\right)^{-1/2}$ (respective vector $d$ is as in (22), $e$ of (37) is such that $e^{(1)} = e^{(2)} = 0$ and condition (24), i.e., inequality $\mathrm{Var}(\delta - \beta\varepsilon) > 0$, is satisfied on account of **(A)**), one concludes convergence to $N(0,1)$ for the initial left hand side in (65). $\qquad\square$

*Proof of the* (b) *part of Theorem* 3. In view of (a) of Theorem 3, the proof reduces to establishing convergence

$$\frac{\sum_{i=1}^{n}(\widetilde{v}_i(j,n) - \overline{\widetilde{v}(j,n)})^2}{\sum_{i=1}^{n}(v_i(j,n) - \overline{v(j,n)})^2} \xrightarrow{P} 1, \quad n \to \infty, \quad j = 1 \text{ and } 2. \tag{69}$$

For $j = 1$ and 2, the expressions for $\widetilde{v}_i(j, n)$ are similar and hence, (69) is shown here for $j = 2$ only.

As $n \to \infty$, from Corollary 1 for $\widehat{\beta}_{2n}$, (60), and WLLN,

$$(\beta - \widehat{\beta}_{2n})^2 S_{xx} = \frac{O_P(1)}{n \ell_\xi^2(n)} \left( \ell_\xi^2(n) O_P(1) + \theta \right) = o_P(1) \quad \text{and} \quad \left( \frac{\overline{x}}{U(2, n)} \right)^2 = \frac{O_P(1)}{\ell_\xi^4(n)}. \tag{70}$$

By (62) and (70),

$$\begin{aligned}
\frac{1}{n} \sum_{i=1}^n & \left( (\widetilde{v}_i(2, n) - \overline{\widetilde{v}_i(2, n)}) - ((v_i(2, n) - \overline{v_i(2, n)})) \right)^2 \\
= \frac{1}{n} \sum_{i=1}^n & \left( (\beta - \widehat{\beta}_{2n})(x_i - \overline{x}) \right. \\
& \left. - \frac{\overline{x}}{U(2, n)} \left( (\widetilde{u}_i(2, n) - \overline{\widetilde{u}(2, n)}) - (u_i(2, n) - \overline{u(2, n)}) \right) \right)^2 \\
\leq 2(\beta - \widehat{\beta}_{2n})^2 S_{xx} & + 2 \left( \frac{\overline{x}}{U(2, n)} \right)^2 \frac{1}{n} \sum_{i=1}^n (\widetilde{u}_i(2, n) - u_i(2, n))^2 \\
= o_P(1) & + \frac{O_P(1)}{\ell_\xi^4(n)} \ell_\xi^2(n) \, o_P(1) = o_P(1),
\end{aligned}$$

that combined with the Cauchy-Schwarz inequality results in (69) for $j = 2$. $\quad\square$

*Proof of Corollary 1 for $\widehat{\alpha}_{jn}$.* Results from the (a) part of Theorem 3, (64), (68) and (57) with $\eta_i(n) = v_i(j, n)$ that are as in the proof of the (a) part of Theorem 3. $\quad\square$

*Proof of Theorem 4.* Below we derive (16) for $k = 1$ only that corresponds to the Studentized CLT in the (a) part of Theorem 1 with $j = 1$ that, according to Theorem 2, also holds true in the SEIVM (1) with an unknown intercept. The proof for $k = 2$ is similar.

Consider the set

$$C_{1n}(\beta) := \left\{ \beta : \frac{\sqrt{n} \, |S_{xy} - \mu| \, |\widehat{\beta}_{1n} - \beta|}{\left( \sum_{i=1}^n \left( (s_{i,yy} - S_{yy}) - \beta(s_{i,xy} - S_{xy}) \right)^2 / (n-1) \right)^{\frac{1}{2}}} \leq z_{\gamma/2} \right\}. \tag{71}$$

On account of the CLT in the (a) part of Theorem 1 with $j = 1$,

$$P\left( C_{1n}(\beta) \right) \to 1 - \gamma, \quad n \to \infty. \tag{72}$$

It is easy to see that

$$C_{1n}(\beta) = \{ \beta : Q_{1n}(\beta, z_{\gamma/2}) \leq 0 \}, \tag{73}$$

with the quadratic function in $\beta$

$$
\begin{aligned}
Q_{1n}(\beta, z_{\gamma/2}) = &\left( n(n-1)(S_{xy} - \mu)^2 - z_{\gamma/2}^2 \sum_{i=1}^{n} (s_{i,xy} - S_{xy})^2 \right) \beta^2 \\
&+ \left( -2n(n-1)(S_{xy}-\mu)^2 \widehat{\beta}_{1n} + 2z_{\gamma/2}^2 \sum_{i=1}^{n} (s_{i,yy} - S_{yy})(s_{i,xy} - S_{xy}) \right) \beta \\
&+ n(n-1)(S_{xy}-\mu)^2 \widehat{\beta}_{1n}^2 - z_{\gamma/2}^2 \sum_{i=1}^{n} (s_{i,yy} - S_{yy})^2,
\end{aligned}
\tag{74}
$$

whose respective discriminant $D_1(n, z_{\gamma/2})$ is as in (19) with $k = 1$. It is crucial to define the sign of $D_1(n, z_{\gamma/2})$ and that of the coefficient of $\beta^2$ in $Q_{1n}(\beta, z_{\gamma/2})$ in order to proceed. On account of $E|\xi|^{8/3} < \infty$ and the Marcinkiewicz-Zygmund law of large numbers, as $n \to \infty$, $\sum_{i=1}^{n} \xi_i^4 / n^{3/2} \xrightarrow{a.s.} 0$ that combined with WLLN leads to $\sum_{i=1}^{n} s_{i,\xi\xi}^2 / n^{3/2} \xrightarrow{a.s.} 0$, and then, to

$$
\sum_{i=1}^{n} s_{i,yy}^2 / n^{3/2} \xrightarrow{a.s.} 0 \quad \text{and} \quad \sum_{i=1}^{n} s_{i,xy}^2 / n^{3/2} \xrightarrow{a.s.} 0.
\tag{75}
$$

On using (75), (61) and Corollary 1 for $\widehat{\beta}_{1n}$, as $n \to \infty$,

$$
\begin{aligned}
&\mathrm{sign}(D_1(n, z_{\gamma/2})) \\
&= \mathrm{sign}\left( (S_{xy} - \mu)^2 \frac{\sum_{i=1}^{n} \left( (s_{i,yy} - S_{yy}) - \widehat{\beta}_{1n}(s_{i,xy} - S_{xy}) \right)^2}{n} \right. \\
&\qquad\qquad - \frac{z_{\gamma/2}^2}{n-1} \left( \frac{\sum_{i=1}^{n} (s_{i,yy} - S_{yy})^2 \sum_{i=1}^{n} (s_{i,xy} - S_{xy})^2}{n^2} \right. \\
&\qquad\qquad\qquad\qquad \left.\left. - \frac{\left( \sum_{i=1}^{n} (s_{i,yy} - S_{yy})(s_{i,xy} - S_{xy}) \right)^2}{n^2} \right) \right) \\
&\geq \mathrm{sign}\left( (S_{xy} - \mu)^2 \frac{\sum_{i=1}^{n} \left( (s_{i,yy} - S_{yy}) - \widehat{\beta}_{1n}(s_{i,xy} - S_{xy}) \right)^2}{n} \right. \\
&\qquad\qquad \left. - z_{\gamma/2}^2 \frac{\sum_{i=1}^{n} s_{i,yy}^2}{n\sqrt{n-1}} \frac{\sum_{i=1}^{n} s_{i,xy}^2}{n\sqrt{n-1}} \right) \xrightarrow{P} 1
\end{aligned}
$$

and

$$
\begin{aligned}
&\mathrm{sign}\left( n(n-1)(S_{xy} - \mu)^2 - z_{\gamma/2}^2 \sum_{i=1}^{n} (s_{i,xy} - S_{xy})^2 \right) \\
&\geq \mathrm{sign}\left( (S_{xy} - \mu)^2 - z_{\gamma/2}^2 \frac{\sum_{i=1}^{n} s_{i,xy}^2}{n(n-1)} \right) \xrightarrow{P} 1.
\end{aligned}
$$

Hence, if

$$C_{2n} := \{\text{sign}(D_1(n, z_{\gamma/2})) = 1\} \ \text{ and}$$

$$C_{3n} := \left\{ \text{sign}\left( n(n-1)(S_{xy} - \mu)^2 - z_{\gamma/2}^2 \sum_{i=1}^{n}(s_{i,xy} - S_{xy})^2 \right) = 1 \right\}, \ (76)$$

then, due to the inequalities

$$P(C_{2n}) \geq P(C_{2n} \cap C_{3n}) \geq P(C_{2n}) + P(C_{3n}) - 1, \tag{77}$$

we have

$$P(C_{2n} \cap C_{3n}) \to 1, \quad n \to \infty. \tag{78}$$

Now, we consider the CI for $\beta$ in (16) with $k = 1$, and denote it by

$$C_n(\beta) := \{\beta : B_1^1(n, z_{\gamma/2}) \leq \beta \leq B_1^2(n, z_{\gamma/2})\}, \tag{79}$$

where $B_1^1(n, z_{\gamma/2})$ and $B_1^2(n, z_{\gamma/2})$ of (17) are the respective smaller and bigger real roots of the quadratic function $Q_{1n}(\beta, z_{\gamma/2})$ of (74) whose discriminant and coefficient of the quadratic term are positive. By (72), (73), (78) and inequalities similar to those in (77),

$$P(C_n(\beta)) = P\Big(C_{1n}(\beta) \cap C_{2n} \cap C_{3n}\Big) + P\Big(C_n(\beta) \cap (C_{2n} \cap C_{3n})^c\Big) \to 1 - \gamma, \quad n \to \infty.$$

$\square$

## Acknowledgements

This paper is based on parts of the author's Ph.D. thesis Martsynyuk (2005), written under the supervision of Miklós Csörgő, and on parts of Martsynyuk (2004). The author is deeply grateful to him for his invaluable guidance throughout her Ph.D. studies, and also for his helpful suggestions and comments that improved the presentation of this exposition.

## References

Carroll, R.J. and Ruppert, D. (1996). The use and misuse of orthogonal regression estimation in linear errors-in-variables models. *American Statistician.* **50** 1-6.

Carroll, R.J., Ruppert, D., Stefanski, L.A. and Crainiceanu, C.M. (2006). *Measurement Error in Nonlinear Models: A Modern Perspective.* 2nd ed. Chapman and Hall/CRC, Boca Raton, London and New York. MR2243417

Cheng, C.-L. and Tsai, C.-L. (1995). Estimating linear measurement error models via M-estimators, in: Mammitzsch, V. and Schneeweiss, H., eds. *Symposia Gaussiana: Proceedings of Second Gauss Symposium, Conference B: Statistical Sciences* (Walter de Gruyter, Berlin) pp. 247–259. MR1355331

Cheng, C.-L. and Van Ness, J.W. (1999). *Statistical Regression with Measurement Error.* Arnold, London. MR1719513

Csörgő, M., Szyszkowicz, B. and Wang, Q. (2004). On weighted approximations and strong limit theorems for self-normalized partial sums processes, in: Horváth, L., Szyszkowicz, B., eds. *Asymptotic Methods in Stochastics: Festschrift for Miklós Csörgő* (Fields Institute Communications **44**, AMS, Providence, Rhode Island) pp. 489–521. MR2108953

Fuller, W.A. (1987). *Measurement error models.* Wiley, New York. MR0898653

Giné, E., Götze, F. and Mason, D.M. (1997). When is the Student $t$−statistic asymptotically standard normal? *Ann. Prob.* **25** 1514-1531. MR1457629

Gleser, L.J. (1987). Confidence intervals for the slope in a linear error-in-variables regression model, in: Gupta, R., ed. *Advances in Multivariate Statistical Analysis* (D. Reidel, Dordrecht) pp. 85–109. MR0920428

Gleser, L.J. and Hwang, J.T. (1987). The nonexistence of $100(1 - \alpha)\%$ confidence sets of finite expected diameter in error-in-variables and related models. *Ann. Statist.* **15** 1351–1362.

Gnedenko, B.V. and Kolmogorov, A.N. (1954). *Limit distributions for sums of independent random variables.* Addison-Wesley, Reading, MA. MR0062975

Maller, R.A. (1981). A theorem on products of random variables, with application to regression. *Austral. J. Statist.* **23** 177–185. MR0636133

Maller, R.A. (1993). Quadratic negligibility and the asymptotic normality of operator normed sums. *J. Multivariate Anal.* **44** 191–219. MR1219203

Martsynyuk, Yu.V. (2004). Invariance principles via Studentization in linear structural error-in-variables models, *Technical Report Series of the Laboratory for Research in Statistics and Probability* **406-October 2004** (Carleton University-University of Ottawa, Ottawa).

Martsynyuk, Yu.V. (2005). *Invariance Principles via Studentization in Linear Structural and Functional Error-in-Variables Models.* Ph.D. dissertation, Carleton University, Ottawa.

O'Brien, G.L. (1980). A limit theorem for sample maxima and heavy branches in Galton-Watson trees. *J. Appl. Probab.* **17** 539–545. MR0568964

Reiersøl, O. (1950). Identifiability of a linear relation between variables which are subject to errors. *Econometrica* **18** 375-389. MR0038054